\documentclass{jcmlatex}

\usepackage{epsfig,amssymb,amsmath,version}
\usepackage{amssymb,version,graphicx,fancybox,mathrsfs}
\usepackage{url,hyperref}
\usepackage{subfigure}
\usepackage{color}
\usepackage{stmaryrd}
\usepackage{multirow}
\usepackage{booktabs,siunitx}
\usepackage{booktabs}
\usepackage{multicol}
\usepackage{booktabs} 
\usepackage{amsmath}
\usepackage{setspace}
\usepackage{array}
\usepackage{placeins}

\usepackage{enumerate}
\usepackage{enumitem}
\usepackage{tabu}                     
\usepackage{multirow}                 
\usepackage{multicol}                 
\usepackage{multirow}                
\usepackage{float}                    
\usepackage{makecell}                 
\usepackage{booktabs}                 
\usepackage{diagbox,bm}
\usepackage{algorithm,float}
\usepackage{algorithmicx}
\usepackage{algpseudocode}
\usepackage{lipsum}
\usepackage{changepage}
\usepackage{graphicx}
\usepackage{booktabs}
\usepackage{multirow}
\usepackage{adjustbox}
\usepackage{afterpage}

\def\JCMvol{xx}
\def\JCMno{x}
\def\JCMyear{202x}
\def\JCMreceived{xxx}
\def\JCMrevised{xxx}
\def\JCMaccepted{xxx}

\floatname{algorithm}{Algorithm 1}

\newtheorem{thm}{\bf Theorem}

\newtheorem{rem}{\bf Remark}
[section]

\definecolor{ligreen}{rgb}{0.0, 0.3, 0.0}

\definecolor{darkblue}{rgb}{0.0, 0.0, 0.55}

\definecolor{anti-flashwhite}{rgb}{0.55, 0.57, 0.68}

\newcommand{\bs}[1]{\boldsymbol{#1}}

\graphicspath{{./Figures/} } 

\makeatletter
\newenvironment{breakablealgorithm}
{
	\begin{center}
		\refstepcounter{algorithm}
		\hrule height.8pt depth0pt \kern2pt
		\renewcommand{\caption}[2][\relax]{
			{\raggedright\textbf{\ALG@name~\thealgorithm} ##2\par}%
			\ifx\relax##1\relax 
			\addcontentsline{loa}{algorithm}{\protect\numberline{\thealgorithm}##2}%
			\else 
			\addcontentsline{loa}{algorithm}{\protect\numberline{\thealgorithm}##1}%
			\fi
			\kern2pt\hrule\kern2pt
		}
	}{
	\kern2pt\hrule\relax
\end{center}
}
\makeatother
\begin{document}

\markboth{Improved Adaptive Orthogonal Basis Method for Multiple Solutions}

\title{An Improved Adaptive Orthogonal Basis Deflation Method for Multiple Solutions with Applications to Nonlinear Elliptic Equations in {\color{black}Varying Domains}}

\author{Yangyi Ye
\thanks{School of Mathematics and Physics, University of South China, Hengyang, China \\ Email: yeyangyi0911@163.com}
\and
 Lin Li \footnotemark[1] 
\thanks{School of Mathematics and Physics, University of South China, Hengyang, China \\ Email:   lilinmath@usc.edu.cn}
Pengcheng Xie \footnotemark[2]
\thanks{Applied Mathematics and Computational Research Division, Lawrence Berkeley National Laboratory, 1 Cyclotron Road, Berkeley, 94720, CA, USA\\ Email: pxie@lbl.gov,pxie98@gmail.com}
Haijun Yu
\thanks{School of Mathematical Sciences, University of Chinese Academy of Sciences, Beijing 100049, China \\
LSEC \& ICMSEC, Academy of Mathematics and Systems Science, Chinese Academy of Sciences, Beijing 100190, China\\
Email: hyu@lsec.cc.ac.cn}}  

\maketitle

\footnotetext[1]{Corresponding authors.}
\footnotetext[2]{The work was done before P. Xie joined LBNL.}
\begin{abstract}
Multiple solutions are common in various non-convex problems arising from industrial and scientific computing. Nonetheless, understanding the nontrivial solutions' qualitative properties seems limited, partially due to the lack of efficient and reliable numerical methods. In this paper, we design a dedicated numerical method to explore these nontrivial solutions further. We first design an improved adaptive orthogonal basis deflation method by combining the adaptive orthogonal basis method with a bisection-deflation algorithm. {\color{black}We then apply the proposed new method to study the impact of domain changes on multiple solutions of certain nonlinear elliptic equations.} {\color{black}When the domain varies from a circular disk to an elliptical disk, the corresponding functional value changes dramatically for some particular solutions}, which indicates that these nontrivial solutions in the circular domain may become unstable in the elliptical domain. 
Moreover, several theoretical results on multiple solutions in the existing literature are verified. For the nonlinear sine-Gordon equation with parameter $\lambda$,  nontrivial solutions are found for $\lambda > \lambda_2$, here $\lambda_2$ is the second eigenvalue of the corresponding linear eigenvalue problem. {\color{black}For the singularly perturbed Ginzburg-Landau equation, highly concentrated solutions are numerically verified, suggesting that their convergent limit is a delta function when the perturbation parameter goes to zero}.
\end{abstract}



\section{Introduction}\label{sect1}

In nonlinear sciences, e.g., physics, mechanics, and biology, the following nonlinear equation is often seen:
\begin{equation}\label{20240110eq1.1}
-\varepsilon^2 \Delta {u}({\bs x}) = f({\bs x}, {u}) \quad\quad\,  \textrm{in}\; \Omega,\\
\end{equation}
where $\varepsilon$ is a real constant, $\Omega$ is a bounded domain in $\mathbb{R}^{d}$ ($d = 1, 2,\text{ or }3$) with a regular boundary $\partial \Omega$. $f$ is a {\color{black}nonlinear function}. Some typical requirements on the regularity and growth of $f({\bs x}, {u})$ are often assumed as follows~(see e.g. \cite{1983Brezis, 1973Rabinowi, Berger1994InfiniteDM, 1991Wang}):
\begin{enumerate}
  \item[${\bs 1})$] $f({\bs x}, {u})$ is locally Lipschitz continuous in $\bar{\Omega}\times \mathbb{R}$;
  \item[${\bs 2})$] For $d\geq 2$, $f({\bs x}, {u}) \leq C_1 + C_2|{u}|^{p}$, where $C_1$ and $C_2$ are constants, $0\leq p \leq \hat{p} = \frac{d+2}{d-2}$;
  \item[${\bs 3})$] 
  There exist a $\mu > 2$ and $M>0$, such that $0 < \mu V({\bs x}, {u}) \leq {u}f({\bs x}, {u})$, for $|{u}|\geq M $, where
  $$V({\bs x}, {u}) = \int_{0}^{{u}}f({\bs x}, v)dv ;$$

  \item[${\bs 4})$] $f({\bs x}, 0) = 0$, or other conditions in a neighborhood of the origin.
\end{enumerate}
When $0<\varepsilon \ll 1$, \eqref{20240110eq1.1} becomes the singularly perturbed semilinear elliptic boundary value problem. 
With a proper boundary condition, e.g. the homogeneous Dirichlet boundary conditions, the corresponding variational function $J: H^{1}_{0}(\Omega)\rightarrow \mathbb{R}$ can be defined by
\begin{equation}\label{20240114eq1.1}
J({u}) = \int_{\Omega}\Bigl{(}\frac{1}{2}\varepsilon^2|\nabla{u}|^{2} - V({\bs x}, {u})\Bigr{)}d{\bs x}.
\end{equation}

Equations of the form \eqref{20240110eq1.1} have attracted significant attention from researchers in different research areas. The existence and multiplicity of solutions of \eqref{20240110eq1.1} have been studied mathematically, see, e.g. \cite{1983Brezis,GUI19991,Chen2003ANO}. Due to the complicated structure of these multiple solutions, computing these solutions often encounters inherent difficulties. Moreover, the choice of the initial guess or basis function plays a crucial role in {\color{black}effectively computing} multiple solutions of \eqref{20240110eq1.1}. To this end, we shall propose an improved adaptive orthogonal basis deflation method for computing multiple solutions, and use it to explore the effect of varying geometry on multiple solutions.  

Before presenting our study, we first give a brief introduction to the research background and related works. Based on the variational functional, i.e. \eqref{20240114eq1.1}, the existence of the solutions with Morse index $1$ was first proved by Ambrosetti and Rabinowitz~\cite{1973Rabinowi} using the mountain pass theorem. Later, Wang~\cite{1991Wang} verified by linking and Morse-type arguments that equation \eqref{20240110eq1.1} with $\varepsilon = 1$ has at least three nontrivial solutions, and the third nontrivial solution is a sign-changing solution with Morse index 2. As a result, some interesting algorithms for computing multiple solutions inspired by these studies have been developed. To be specific, in 1993, based on the mountain pass theorem, a mountain pass algorithm (MPA) ~\cite{choi1993mountain} was designed to compute multiple solutions of \eqref{20240110eq1.1}.
It is worth pointing out that the MPA is only limited to finding two solutions of mountain pass type with Morse index 1 or 0. As an extension of MPA, Ding et al.~\cite{ding1999high}  proposed a high-linking algorithm (HLA) to find sign-changing solutions of \eqref{20240110eq1.1}. Later, Zhou et al. \cite{li2001minimax} proposed a min-max algorithm (MNA) for finding multiple critical points of \eqref{20240114eq1.1}, where some local min-max theorems were established and a two-level local optimization on a solution sub-manifold
is suggested. MNA inspired several later works, e.g., \cite{2005A, Zhou2007Nu, 2005Instability, 2011A,xie100810411,Li2018ALM,Li2018,LI2013407,xie2022}. {\color{black}For non-variational type of differential equations, several multiple-solution algorithms have also been proposed}, where the main idea is to use some iteration methods to find multiple solutions of a nonlinear algebraic system (NLAS) from the discretization of the differential equations. The search-extension method (SEM)~\cite{chen2004search} was proposed to compute \eqref{20240110eq1.1} with $\varepsilon= 1$, where the following {\color{black}eigenvalue problems} are first solved
\begin{equation*}
-\Delta \psi_j = \lambda_j \psi_j   \quad  \textrm{in}\;\Omega; 
\qquad  \psi_j=0 \text{ or } \partial_n \psi_j = 0, 
\quad\textrm{on }\partial\Omega. \vspace{-0.1cm}
\end{equation*}
Then, the solution of \eqref{20240110eq1.1} can be approximated by
\begin{equation*}
{u}_{N}({\bs x}) = \sum_{j = 1}^{N}a_j\psi_j({\bs x}) \in S_N,   \quad\quad\quad  S_N := \text{span}\{ \psi_1, \psi_2, \ldots, \psi_N\}.
\end{equation*}
Some improvements on the SEM are given in \cite{xie2006improved}. {\color{black} Later, Chen at el.~\cite{2008ChenZhou} and Allgower et al.~\cite{allgower2006solution, allgower2009application} proposed a homotopy continuation method, respectively,} where a homotopy tracking number $\tau$ was introduced, and the corresponding homotopy function was solved by Newtonian iterations. Several further works have also been developed, e.g., \cite{hao2014bootstrapping, 2018Two,Hao2020AHM,Hao2020SpatialPF,Zhao2022BifurcationAR}. Another important approach is presented by Farrell et al.~\cite{farrell2015deflation}, where the deflation technique and the Newtonian iterations were integrated to solve the NLAS from the discretization of nonlinear partial differential equations. Later, Farrell et al. extended the deflation technique
\cite{Farrell2020,Farrell2017Siam,Farrell2018,Farrell2019}, where the deflation approach was interwoven with nested iteration, creating an efficient method for multiple solutions. It is worth pointing out that the methods mentioned above are based on the Newtonian iteration.
As we know, the Newtonian iteration is sensitive to the initial guess and requires that the Jacobian matrix has full column rank, making it not numerically robust, nor computationally efficient.
To improve it, an efficient spectral trust-region method \cite{2022Two} is proposed, where the classical trust-region method is introduced to replace the Newtonian iteration.
Numerical results illustrate that the classical trust region method has noteworthy advantages compared with the Newtonian method. So, the trust region method will be adopted again in the current study. Recently, by considering the equation's nature and structural characteristics of the solution, an efficient adaptive orthogonal basis method for computing multiple solutions of \eqref{20240110eq1.1} with polynomial nonlinearities is first developed~\cite{Li2024}. Different from conventional practices of predefining candidate basis pools, a set of adaptive orthogonal basis functions was constructed. Meanwhile, the companion matrix technique is used to generate good initial guesses for computing multiple solutions. Unfortunately, the adaptive orthogonal basis method is confined to systems with only polynomial nonlinearities, which greatly limits its application to more generalized systems. Therefore, one of the major aims of this paper is to extend the adaptive orthogonal basis method to \eqref{20240110eq1.1} with more general nonlinearities, where  
a novel bisection method interwoven with the deflation technique is introduced to provide good initial guesses for computing multiple solutions of \eqref{20240110eq1.1}. {\color{black}As a result, an improved adaptive orthogonal basis deflation method (IAOBDM) for computing multiple solutions of \eqref{20240110eq1.1} is proposed}.

Multiple solutions of nonlinear elliptic equations have many practical applications~(see e.g.\cite{GUI19991,Cao1999OnTE,GUI200047}), including topology optimizations in structural and additive manufacturing~(see e.g. \cite{Farrell2021Siam}). {\color{black}Exploring the impact of domain changes on multiple solutions of} \eqref{20240110eq1.1} is helpful for guiding structural and additive manufacturing. However, {\color{black}very few published literature exist on this topic}. {\color{black}In this paper, we study the influence of domain changes on multiple solutions of \eqref{20240110eq1.1} by using the IAOBDM.}
In addition, we shall numerically verify some interesting theoretical results on multiple solutions presented in the published literature using the proposed method. For example, for the singularly perturbed boundary value problem \eqref{20240110eq1.1} with the homogeneous Neumann boundary condition, it is known that it has a solution ${u}$ with exactly one single peak (i.e., the local maximum of ${u}$). Moreover, if the single peak is on the boundary of $\Omega$, {\color{black}then it is located at the point where the mean curvature of the boundary $\partial\Omega$ reaches its maximum value.} These results have not been numerically verified, which limits the further application of multiple solution analysis.

The remainder of this paper is organized as follows. In Section \ref{sect2}, we describe the spectral Legendre--Fourier scheme used to discretize the equation \eqref{20240110eq1.1} defined in an elliptic geometry. The IAOBDM is designed and presented in section \ref{sect3}. In Section \ref{sect4}, ample numerical experiments are carried out to demonstrate the efficiency of this method and to show the effect of varying geometry $\Omega$ on multiple solutions of \eqref{20240110eq1.1}. Finally, we end the paper with some remarks in Section \ref{sect5}.\vspace{-0.7cm}

\section{A Legendre--Fourier scheme for elliptic equations in an elliptic domain}\label{sect2}
To numerically study the effect of varying geometry $\Omega$ on multiple solutions of \eqref{20240110eq1.1}, the first step is to provide an efficient discretization scheme. Here, we adopt a Legendre--Fourier scheme for \eqref{20240110eq1.1} in an elliptic domain. Let
\begin{equation}
\Omega = \bigl\{ (x, y): \tfrac{x^2}{a^2} + \tfrac{y^2}{b^2} \leq 1 \bigr\}.
\end{equation}
The weak formulation of \eqref{20240110eq1.1} with homogeneous Dirichlet boundary condition and $\varepsilon = 1$ is to find $u\in H^1_{0}(\Omega)$ such that
\begin{equation}\label{20240101eq2.1}
\int_{\Omega}\nabla {u} \nabla {v} dxdy  =  \int_{\Omega} f({\bs x}, {u}){v} dxdy,  \quad   \forall {v} \in H^1_0(\Omega), 
\end{equation}
where ${\bs x} = (x, y)$.
We use polar transformation $x = a r\cos\theta, y = b r \sin\theta$ to transform the Dirichlet problem into polar coordinate form:
\begin{equation}\label{eq15}
\begin{cases}\vspace{0.15cm}
\omega_1 u_{rr} + \dfrac{1}{r}( \omega_2 u_{r} - \omega_3 u_{r\theta})
+\dfrac{1}{r^2}(\omega_3 u_{\theta} + \omega_2 u_{\theta\theta}) + f(r, \theta, u) = 0, \quad  \textrm{in}\; \Omega = (0, 1)\times [0, 2 \pi),\\
u(1, \theta) = 0, \quad\quad\quad   u \textrm{ periodic in } \theta,
\end{cases}
\end{equation}
where
\begin{equation}
\omega_1 = \frac{\cos^2\theta}{a^2} + \frac{\sin^2\theta}{b^2},  \qquad
\omega_2 = \frac{\cos^2\theta}{b^2} + \frac{\sin^2\theta}{a^2},  \qquad
\omega_3 = \sin2\theta(\frac{1}{a^2}-\frac{1}{b^2}).
\end{equation}
The corresponding weak formulation \eqref{20240101eq2.1} becomes
\begin{equation}\label{20240101eq2.5}
\begin{split}
\mathcal{L}(u, v) := & \int_{\Omega}(\omega_{3}u_{\theta}v_{r}-\omega_2 u_{r}v-\omega_1(u_{r}v+ru_{r}v_{r}))drd\theta\\
~ & + \int_{\Omega}\frac{1}{r}u_{\theta}(\omega_3 v - \omega_2 v_{\theta})drd\theta = \int_{\Omega}fvdrd\theta, \quad  \forall v \in H^1_0(\Omega).
\end{split}
\end{equation}
Obviously, from \eqref{eq15} or \eqref{20240101eq2.5}, we can find that a singularity point appears at $r = 0$, indicating that additional polar conditions should be proposed to obtain the desired solution regularity at $r=0$. 
In the polar coordinate system, 
the solution $u$ to \eqref{20240101eq2.5} is expanded as
\begin{equation}
u(r, \theta) = \sum_{|k| = 0}^{\infty}\bigl(u_{1k}(r)\cos k\theta + u_{2k}(r)\sin k\theta\bigr).
\end{equation}
For solution $u\in C^\infty(\Omega)$, the coefficient functions $u_{1 k}(r)$ or $u_{2 k}(r)$ should satisfy the following polar conditions (see \cite{1983Orszag}):
\begin{equation}\label{20240101eq2.3}
u_{1k}(r) = O(r^{|k|}),\quad\quad  u_{1k}(r) = O(r^{|k|}),\quad\quad  \textrm{as}\; r\rightarrow 0, \quad |k|=1, 2, \cdots.
\end{equation}
However, it is not easy to impose \eqref{20240101eq2.3} in a numerical scheme. 
In fact, for $u \in H_0^1(\Omega)$ as in \eqref{20240101eq2.5}, only the following polar conditions are necessary (see e.g. \cite{1997shenjie, shen_approximation_2012}):
\begin{equation}\label{20240101eq2.4}
u_{1k}(0) = u_{2k}(0) = 0   \quad\quad\quad   \textrm{for} \;\; k \neq 0,
\end{equation}
which is termed as the essential polar condition for \eqref{20240101eq2.5}, while polar conditions in \eqref{20240101eq2.3} is termed as natural or nonessential. In this paper, we use \eqref{20240101eq2.4}, thus $u(r, \theta)$ should satisfy the following boundary conditions:
\begin{equation}\label{20240111eq2.1}
\partial_{\theta} u(0, \theta) = 0, \quad\quad   u(1, \theta) = 0.
\end{equation}
A Fourier spectral approximation  to \eqref{20240101eq2.5} can be defined by
\begin{equation}\label{20240111eq2.2}
u_{M, N}(r, \theta) = \sum_{i=1}^{M}\sum_{j=0}^{N-2}(\alpha_{ij}\sin(i\theta)\hat{\phi}_{j}(r)+ \beta_{ij}\cos(i\theta)\hat{\phi}_{j}(r)) + \sum_{l=0}^{N-1}\gamma_{l}\hat{\varphi}_{l}(r),
\end{equation}
where $\alpha_{ij}, \beta_{ij}$ and $\gamma_{l}$ are unknown variables to be determined, and $\hat{\phi}_{j}(r)$ and $\hat{\varphi}_{l}(r)$ satisfy
\begin{equation}\label{20240111eq2.3}
\hat{\phi}_{j}(0) = \hat{\phi}_{j}(1) = 0, \quad
0\leq j\leq N-2, 
\qquad 
\hat{\varphi}_{l}(1) = 0,  \quad 0\leq l \leq N-1.
\end{equation}
Introducing $t = 2r-1$, we convert $r\in[0,1]$ into $t\in [-1,1]$.
The corresponding problem in $(t, \theta)$ coordinates reads:
\begin{equation}\label{20240111eq2.4}
\begin{cases}\vspace{0.15cm}
\omega_1 {u}_{tt} + \dfrac{1}{2(t+1)}( \omega_2 {u}_{t} - \omega_3 u_{t\theta})
+\dfrac{1}{(t+1)^2}(\omega_3 {u}_{\theta} + \omega_2 {u}_{\theta\theta}) + \frac{1}{4}f({u}) = 0,\; \textrm{in}\; \Omega = (-1, 1)\times [0, 2\pi),\\
u(1, \theta) = 0, \quad
\partial_\theta u(-1, \theta)=0, \quad   
u \textrm{ periodic in }\theta.
\end{cases}
\end{equation}
In the Legendre--Fourier spectral method, we use the following expansion to approximate the solution
\begin{equation}\label{20240112eq2.1}
\hat{u}_{M, N}(t, \theta) = \sum_{i=1}^{M}\sum_{j=0}^{N-2}(\alpha_{ij}\sin(i\theta)\phi_{j}(t)+ \beta_{ij}\cos(i\theta)\phi_{j}(t)) + \sum_{l=0}^{N-1}\gamma_{l}\varphi_{l}(t),
\end{equation}
The $\phi_{j}(t)$ and $\varphi_{l}(t)$ in \eqref{20240112eq2.1} are taken as
\begin{equation}
\phi_{j}(t) = c_{j}(L_{j}(t) - L_{j+2}(t)),\quad\;   \varphi_{l}(t) = L_{l}(t) - L_{l+1}(t),
\end{equation}\vspace{-0.25cm}
where $L_j(t)$ stands for Legendre polynomial of degree $j$, and 
\begin{equation}
c_{j} = \frac{1}{\sqrt{4j+6}},\quad\quad  \int_{-1}^{1}L_{m}(t)L_{n}(t)dt = \frac{1}{2n+1}\delta_{mn}.
\end{equation}
Denoting
\begin{equation}
{\bs \xi} = (\alpha_{10}, \ldots, \alpha_{M(N-2)}, \beta_{10}, \ldots, \beta_{M(N-2)}, \gamma_0, \ldots, \gamma_{N-1})^{\top},
\end{equation}
the Galerkin approximation of \eqref{20240111eq2.4} is to find ${\bs \xi}$ such that
\begin{equation}\label{20240112eq2.2}
\mathcal{F}(\hat{u}_{M, N}, \psi) = 0,   \quad\quad  \forall \psi \in X_{M, N},
\end{equation}
where $X_{M, N} := \text{span}\bigl\{\sin(i\theta)\phi_{j}(r),\; \cos(i\theta)\phi_{j}(r),\; \varphi_{l}(r), \; 1\leq i\leq M,\; 0\leq j \leq N-2,\; 0\leq l\leq N-1\; \bigr\}$. $\mathcal{F}$ is a form similar to \eqref{20240101eq2.5}.
Let
\begin{equation*}
a_{ij} = \int_{-1}^{1}(t+1)\phi'_{i}(t)\phi'_{j}(t)dt,  \quad\;   b_{ij} = \int_{-1}^{1}\frac{1}{t+1}\phi_{i}(t)\phi_{j}(t)dt, \quad\;  c_{ij} = \int_{-1}^{1}(t+1)\phi_{i}(t)\phi_{j}(t)dt,
\end{equation*}
\begin{equation*}
d_{ij} = \int_{-1}^{1}(t+1)\varphi'_{i}(t)\varphi'_{j}(t)dt,\quad\quad\quad
e_{ij} = \int_{-1}^{1}(t+1)\varphi_{i}(t)\varphi_{j}(t)dt.
\end{equation*}
They are the $t$-direction matrices related to the linear part of the Galerkin discretized equation. It is easy to show that the above matrices are symmetric and sparse, their nonzeros in the up-right parts are given as
\begin{equation*}
a_{ij} = \begin{cases}
2i+4, & j = i+1\\
4i+6, & j = i,
\end{cases}
\quad\;
b_{ij} = \begin{cases}\vspace{0.15cm}
-\frac{2}{i+2},  & j = i+1,\\
\frac{2(2i+3)}{(i+1)(i+2)}, & j = i,
\end{cases}\quad\;
d_{ij} = \begin{cases}\vspace{0.15cm}
2i+2,   & j = i,\\
0,   &  \textrm{otherwise},
\end{cases}
\end{equation*}
\begin{equation*}
c_{ij} = \begin{cases}\vspace{0.15cm}
-4\frac{i+3}{2(2i+5)(2i+7)},   & j = i+3,\\\vspace{0.15cm}
-\frac{2}{2i+5},   &  j = i+2,\\\vspace{0.15cm}
\frac{2}{(2i+1)(2i+5)} + \frac{2(i+3)}{(2i+5)(2i+7)},   & j = i+1,\\\vspace{0.15cm}
\frac{2}{2i+1} + \frac{2}{2i+5},  & j = i,
\end{cases}\quad\quad
e_{ij} = \begin{cases}\vspace{0.15cm}
-\frac{2(i+2)}{(2i+3)(2i+5)},  & j = i+2,\\\vspace{0.15cm}
-\frac{4}{(2i+1)(2i+3)(2i+5)},  & j = i+1,\\\vspace{0.15cm}
\frac{4(i+1)}{(2i+1)(2i+3)},   & j = i.
\end{cases}
\end{equation*}
Equation \eqref{20240112eq2.2} can be written in a NLAS as follow:\vspace{-0.2cm}
\begin{equation}\label{eq2.2.4.1}
{\bs F} =\big({F}_{1}, \; {F}_{2}, \; \cdots, \; {F}_{n}\big)^{\top} = {\bs 0},
\end{equation}
where $n = 2M(N-1) + N$ and $\{F_i\}_{i=1}^{n}$ are algebraic functions of the unknown vector ${\bs \xi}$.

\section{An improved adaptive orthogonal basis deflation method}\label{sect3}
In this section, we develop the IAOBDM method to calculate multiple solutions of \eqref{eq2.2.4.1} (or \eqref{20240110eq1.1}) by extending the adaptive orthogonal basis method proposed in \cite{Li2024}. To be specific, when $f({\bs x}, {u})$ in \eqref{20240110eq1.1} is not a polynomial of ${u}$, the adaptive orthogonal basis method proposed in \cite{Li2024} is not applicable due to the limitation of the companion matrix technique. The main difficulty is that {\color{black}it is hard to find all multiple roots of a nonlinear algebraic system}. To our knowledge, so far none of the existing algorithms can work well. Here we design a robust bisection-deflation algorithm for finding multiple roots of a nonlinear algebraic equation. 

We also give a brief introduction to 
the trust region method adopted to find solutions to the NLAS.  The IAOBDM is a combination of the bisection-deflation algorithm, trust region method and adaptive orthogonal basis method.

\subsection{A bisection-deflation algorithm for multiple roots of scalr algebraic equations}\label{sect3.1} We assume that a scalar nonlinear algebraic equation is given as
\begin{equation}\label{2024080801}
\omega(x) = 0,    \quad\quad\quad   x\in \mathbb{R},
\end{equation}
where $\omega$ is a continuous function.
In general, due to the limitation of the classical Newtonian method, it is very difficult to find {\color{black}all multiple roots of \eqref{2024080801}}. 
{\color{black}When combing the deflation and the trust region method presented in \cite{2022Two}}, the corresponding nonlinear iteration will become more and more complicated as the number of multiple solutions increases.
Here, we apply a bisection method presented in \cite{1969DekkerTJ} to reduce the computational complexity. We interwoven the bisection method with the deflation technique to compute multiple roots of \eqref{2024080801}. The detailed bisection-deflation method is described in \textbf{Algorithm} \ref{2024080901}, in which the following pre-defined functions are used.
\begin{equation}\label{eq:fun-omegas}
    \begin{aligned}
        \omega_1 (b, a) &= \begin{cases}
            b- \frac{b-a}{\omega(b)-\omega(a)}\omega(b),  & \textrm{if } a\neq b \text{ and } \omega(b) \neq \omega(a),\\
            \infty,  & \textrm{if } a \neq b \text{ and } \omega(b) = \omega(a) \neq 0,\\
            b,  & \textrm{otherwise},
    \end{cases}\\
        \omega_{2}(b, c) &= b + \textrm{sign}(c-b)\times \textbf{eps}(b), \\
        \omega_{3}(b, c) &= \frac{1}{2}(b+c),  \\    
        \omega_{4}(x, b, c) &= \begin{cases}
            x,   & \textrm{ if }\, x\, \textrm{ is }\, \textrm{between }\,   \omega_{2}(b,c)\,\textrm{ and }\, \omega_{3}(b,c),\\
            \omega_{2}(b,c),  & \textrm{ if }\; |x - b|\leq \textbf{eps}(b),\\
            \omega_{3}(b,c),  & \textrm{ otherwise},
    \end{cases}
    \end{aligned}
\end{equation}
where $\textbf{eps}(\cdot)$ represents the relative accuracy function in the floating-point number system used.

\begin{center}
\begin{breakablealgorithm}
\setstretch{1.2}
\caption{An efficient bisection-deflation algorithm for computing roots of \eqref{2024080801}\label{2024080901}}
\begin{algorithmic}[0]
\State {\bf Step 1.}
Initialize the solution set as $\mathcal{S}=\emptyset$. 

\State {\bf Step 2.}\vspace{-0.5cm}
\begin{adjustwidth}{1.5cm}{0.5cm}
(\textbf{Initialization a new root search}).
Try to find two distinct values $x_{0}$ and $x_{1}$ satisfying $\omega(x_{0})\times \omega(x_{1}) \leq 0$ by
a grid search in a region that determined by the regularity and growth of $f$. 

If not succeed, then return $\mathcal{S}$.

If succeeds, then set
\begin{equation}
    \begin{cases}
        b_1 = x_1, \quad a_1 = c_1 = x_0, & \textrm{if }    \omega(x_{1}) \leq \omega(x_{0}),\\
        b_1 = x_0, \quad a_1 = c_1 = x_1, & \textrm{otherwise}.
    \end{cases}
\end{equation}

\end{adjustwidth}

\State {\bf Step 3.}\vspace{-1cm}
\State
\begin{adjustwidth}{1.5cm}{0.5cm}
(\textbf{Iteration}, \textbf{$k = 2, \ldots, $}). 
 The iteration point $x_{k}$ is obtained by $x_{k} = \omega_{4}(\lambda_{k}, b_{k-1}, c_{k-1})$, where $\lambda_{k} = \omega_{1}(b_{k-1}, a_{k-1})$. Let $j$ be the largest integer satisfying $j < k$ and $\omega(x_{j})\times \omega(x_{k}) \leq 0$. Then, $b_{k}, c_{k}$ and $a_{k}$ are defined by
    
\begin{equation}
\begin{cases}
    b_{k} = x_{k}, \quad  c_{k} = x_{j}, \quad  a_{k} = b_{k-1},&  \quad \textrm{ if } \omega(x_{k}) \leq \omega(x_{j}), \\
    b_{k} = x_{j},  \quad a_{k} = c_{k} = x_{k}, & \quad \textrm{ otherwise}.
\end{cases}
\end{equation}
\end{adjustwidth}
\State
\begin{adjustwidth}{1.5cm}{0.5cm}
(\textbf{Check convergence}). If $|b_{k} - c_{k}|\leq 2\, \textbf{\textrm{eps}}(b_{k})$, then 
        set
	\begin{equation}
	x^{*} = b_{n}\; \textrm{or}\; c_{n},
	\end{equation}
	and let $\mathcal{S} = \mathcal{S} \cup \{ x^* \}$, exit the loop.
\end{adjustwidth}
\State {\bf Step 4.}\vspace{-1cm}
\State
\begin{adjustwidth}{1.5cm}{0.5cm}
(\textbf{Deflation}). Use the deflation function $\psi(x, x^{*}) = \frac{1}{|x - x^{*}|^{2}} + 1$ to deflate the known root $x^{*}$, i.e., set $\omega(x)\leftarrow \psi(x, x^{*})\omega(x)$, and then go to \textbf{step 1}.
\end{adjustwidth}
\end{algorithmic}
\end{breakablealgorithm}
\end{center}

\begin{rem} In \textbf{\color{black}Step 4}, the deflation function $\psi$ is introduced {\color{black}to keep the iteration away from a known root so that a new root can be obtained.} 
The reader is referred to \cite{1971Deflation, farrell2015deflation, 2022Two} for more details.
For the asymptotic order of convergence of the bisection part of \textbf{Algorithm} \ref{2024080901}, please see  \cite{1969DekkerTJ} for more details. 
\end{rem}

\subsection{Trust region method for NLAS}

\begin{figure}[h]
\begin{center}	\includegraphics[width=0.28\textwidth]{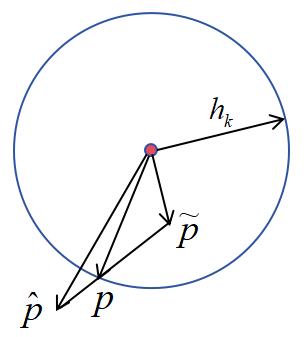}
\caption{Dogleg approximation.}\label{2024082701}
\end{center}
\end{figure}

According to results from previous research \cite{2022Two}, the trust region method has remarkable advantages compared to the classical Newtonian iteration. Therefore, the trust region method is used here to solve \eqref{eq2.2.4.1}. To apply this method, we first rewrite \eqref{eq2.2.4.1} into a minimization problem:
\begin{equation}\label{eq2.2.4.2}
\min_{{\bs \xi} \in {\mathbb R}^{n}}Q({\bs \xi}),\quad Q({\bs \xi}):= \frac{1}{2}\big\|\bs{F}({\bs \xi})\big\|^{2}_{2} = \frac{1}{2}\sum_{i = 1}^{n}{F}^2_{i}({\bs \xi}).
\end{equation}
Our aim is to find some ${\bs \xi}$ to satisfy ${\bs F}({\bs \xi}) \equiv {\bs 0}$, i.e., $Q({\bs \xi}) \equiv 0$. In other words, solve \eqref{eq2.2.4.2} instead of \eqref{eq2.2.4.1} using the trust region method. The main idea is to define a trustable region around the current iteration point. In the $k$-th iteration, a quadratic model to be a suitable representation of the objective function $Q({\bs \xi})$ in \eqref{eq2.2.4.2} is constructed using the Taylor formula:
\begin{equation}\label{eq2.2.4.4}
\begin{split}
& \min_{\bs s\in {\mathbb B}_{h_k}} q^{(k)}({\bs s}) := Q({\bs \xi}_{k}) + {\bs g}({\bs \xi}_{k})^{\top}{\bs s} + \frac{1}{2}{\bs s}^{\top}{\bs G}({\bs \xi}_{k}){\bs s},   \quad\quad  k \geq 0,
\end{split}
\end{equation}
where the trust region ${\mathbb B}_{h_k}:=\{\bs s\in {\mathbb R}^n\,:\,\|\bs s\|\le h_k\}$ and\vspace{0.1cm}
\begin{subequations}\label{New4.4}
\begin{equation}
\bs J({\bs \xi}_{k}) = {\bs F}'({\bs \xi}_{k}) = (\nabla {F}_{1}({\bs \xi}_{k}), \nabla {F}_{2}({\bs \xi}_{k}), \cdots, \nabla F_{n}({\bs \xi}_{k}))^{\top},\vspace{0.1cm}
\end{equation}
\begin{equation}\label{2024082703}
{\bs g}_{k} := {\bs g}({\bs \xi}_{k}) = \nabla Q({\bs \xi}_{k}) = \bs{J}^{\top}({\bs \xi}_{k}) {\bs F({\bs \xi}_{k})},\vspace{0.1cm}
\end{equation}
\begin{equation}\label{2024082704}
{\bs G}_{k} := {\bs G}({\bs \xi}_{k}) = \nabla^2 Q({\bs \xi}_{k}) = {\bs J}^{\top}({\bs \xi}_{k}) \bs J({\bs \xi}_{k}) + \sum_{i = 1}^{n}F_{i}({\bs \xi}_{k})\nabla^2 F_{i}({\bs \xi}_{k}).
\end{equation}
\end{subequations}
Then, the classical dogleg method presented in \cite{Li2024} is used to solve \eqref{eq2.2.4.4} (see Fig. \ref{2024082701}), i.e., the solution ${\bs s}$ to \eqref{eq2.2.4.4} becomes

\begin{equation}\label{2024082702}
{\bs s}_{k}=\begin{cases}\vspace{0.2cm}
-{\bs G}^{-1}_{k}{\bs g}_{k}, &  \textrm{if}\;\; \|\hat{\bs p}\|_{2} \leq h_{k},    \\\vspace{0.2cm}
  -\frac{{\bs g}_{k}}{\|{\bs g}_{k}\|_2} h_{k}, &  \textrm{if}\;\; \|\tilde{\bs p}\|_{2} > h_{k} \;\, \textrm{or}\;\,  \|\hat{\bs p}\|_{2} > h_{k}, \\\vspace{0.2cm}
 {\bs p} := \tilde{\bs p} + \lambda(\hat{\bs p} - \tilde{\bs p}), & \textrm{if}\;\;  \|\hat{\bs p}\|_{2} > h_{k} \;\;
 \textrm{and}\;\;  \|\tilde{\bs p}\|_{2} \leq h_{k},
\end{cases}
\end{equation}
where
\begin{subequations}
\begin{equation*}
\hat{\bs p} = -{\bs G}^{-1}_{k}{\bs g}_{k}  \quad\quad\quad\quad\quad \tilde{\bs p} = -\frac{{\bs g}^{\top}_{k}{\bs g}_k}{{\bs g}^{\top}_{k}{\bs G}_k{\bs g}_k}{\bs g}_k, \vspace{0.2cm}
\end{equation*}
\begin{equation*}
\lambda = \frac{\tilde{\bs p}^{\top}(\hat{\bs p}-\tilde{\bs p})-\sqrt{ (\tilde{\bs p}^{\top}(\hat{\bs p}-\tilde{\bs p}))^2 - \|\hat{\bs p} - \tilde{\bs p}\|^2_2 (\|\hat{\bs p}-\tilde{\bs p}\|^2_2 - h^2_k) }}{\|\hat{\bs p} - \tilde{\bs p}\|^2_2}.
\end{equation*}
\end{subequations}

\vspace{-0.5cm}
When the objective function $Q({\bs \xi})$ is complicated, it may not be easy to compute \eqref{2024082703}-\eqref{2024082704} exactly. Finite differences can be used to calculate approximate derivatives as follows:
\begin{subequations}
\begin{equation}
({\bs g}_k)_{i} \approx \frac{Q({\bs \xi}_k+\epsilon{\bs e}_i) - Q({\bs \xi}_k-\epsilon{\bs e}_i)}{2\epsilon},
\end{equation}
\begin{equation}
({\bs G}_k)_{ij} \approx \frac{Q({\bs \xi}_k+\epsilon{\bs e}_i +\epsilon{\bs e}_j) - Q({\bs \xi}_k+\epsilon{\bs e}_i) - Q({\bs \xi}_k+\epsilon{\bs e}_j) + Q({\bs \xi}_k)}{\epsilon^2},
\end{equation}
\end{subequations}
where $0 < \epsilon \ll 1$ and ${\bs e}_i$ is the $i$-th unit vector. The gradients and Hessians can also be approximated by interpolation \cite{xie2023linesearch,Xie2023LeastN,xie2023dfoto,xie2023twodimensional,
Xie2024}. To measure a good agreement between $q^{(k)}(s_{k})$ and the objective function value $Q({ \bs \xi}_{k} + {\bs s}_{k})$, we should choose a suitable trust region radius $h_{k}$ adaptively, for which a ratio $r_{k}$ should be introduced, and defined by
\begin{equation}\label{rkval}
r_{k} = \frac{Q({\bs \xi}_{k}) - Q({\bs \xi}_{k} + {\bs s}_{k})}{q^{(k)}(\bs 0) - q^{(k)}({\bs s}_{k})}.
\end{equation}
The ratio $r_{k}$ is an indicator of the effectiveness of each step and plays a key role in expanding and contracting $h_{k}$. Compared with the Newtonian method and the line search method, the direction and length of the step are considered simultaneously in the trust region method. The trust region radius $h_{k}$ will be reduced when step ${\bs s}_{k}$ is deemed unacceptable, and a new minimizer should be sought. If trust region radius $h_{k}$ is too large, then the minimizer of the constructed quadratic model may be far from the one of the true objective function $Q({\bs \xi}_{k} + {\bs s}_{k})$ given in \eqref{eq2.2.4.2}. On the contrary, if $h_k$ is too small, a substantial step that can bring the solution much closer to the minimizer of the objective function may be missed. In practical computations, the trust region radius is often adjusted based on the numerical performance in previous iterations. According to \cite{211Wright}, if the constructed quadratic model \eqref{eq2.2.4.4} provides reliable results, the trust region radius $h_{k}$ should be kept or increased to allow larger and more ambitious steps. Otherwise, the trust-region radius $h_{k}$ should be decreased. At each step, the trust region radius $h_{k}$ needs to be reconsidered and changed. The detailed algorithm of the trust region method for solving \eqref{eq2.2.4.2} is summarized in \textbf{Algorithm} \ref{2024080902}.

\begin{breakablealgorithm}
\begin{center}
\setstretch{1.2}
\caption{A trust region method for solving \eqref{eq2.2.4.2}}\label{2024080902}
\begin{algorithmic}[1]
\Require Given ${\bs \xi}_0$, $\epsilon_g > 0$, $\epsilon_v > 0$,  $0<\delta_1 < \delta_2 < 1, 0 < \tau_1 < 1 < \tau_2$ and $h_{0} = \|\bs g_{0}\|$.
\Ensure ${\bs \xi}$
\State \textbf{For} $k = 0,~ 1,~ 2,~ \cdots$
\State \hspace{0.4cm}Compute ${\bs g}_{k}$ and ${\bs G}_{k}$;
\State \hspace{0.4cm}If $\|{\bs g}_{k}\| \leq \epsilon_g$ and $|Q({\bs \xi}_{k})| \leq \epsilon_v$, exit the for-loop and go to line 9;
\State \hspace{0.4cm}Approximately solve the quadratic model (\ref{eq2.2.4.4}) with \eqref{2024082702};
\State
\vspace{-0.5cm}
\begin{adjustwidth}{0.45cm}{0.5cm}
Compute $Q({\bs \xi}_{k} + {\bs s}_{k})$ and $r_{k}$. If $r_{k}\geq \delta _1,$ then ${\bs \xi}_{k+1} = {\bs \xi}_{k} + {\bs s}_{k}$; Otherwise, set ${\bs \xi}_{k+1} = {\bs \xi}_{k}$;
\end{adjustwidth}
\State \hspace{0.4cm}If $r_{k} < \delta_1$, then $h_{k+1} = \tau_1 h_{k}$;
\State \hspace{0.4cm}If $r_{k} > \delta_2$ and $\|{\bs s}_{k}\|$ = $h_{k}$, then $h_{k+1} = \tau_2 h_{k}$; Otherwise, set $h_{k+1} = h_{k}$;
\State \textbf{End}
\State Return ${\bs \xi} = {\bs \xi}_{k}$.
\end{algorithmic}
\end{center}
\end{breakablealgorithm}
In the current paper, we choose $\epsilon_g = \epsilon_v = 10^{-13}, \delta_1 = 0.25, \delta_2 = 0.75, \tau_1 = 0.5,$ and $\tau_2 = 2$. 

It is known that the trust region method given above enjoys desirable global convergence with a local superlinear rate of convergence, as stated in the following theorem.

\begin{thm}[see \cite{sun2006optimization}]\label{thm2024082705} Assume that the objective function $Q({\bs \xi})$ is Lipschitz continuously differentiable on the level set $H := \{{\bs \xi}\in \mathbb{R}^{n}\; :\; Q({\bs \xi}) \leq Q({\bs \xi}_{0})\}, \quad \forall\, {\bs \xi}_{0}\in \mathbb{R}^n$, and the Hessian matrices ${\bs G}_{k}$ are uniformly bounded in 2-norm. Then, $\lim_{k \to +\infty}\inf\|{\bs g}_{k}\| = 0.$ Moreover, if $\lim_{k \to +\infty}{\bs \xi}_k \rightarrow {\bs \xi}^*$ with ${\bs g}({\bs \xi}^{*}) = {\bs 0}$, and ${\bs G}({\bs \xi}^{*})$ is positive definite, the convergence rate of the trust region method is quadratic.
\end{thm}

\subsection{\color{black}Algorithm of IAOBDM}\label{sect3.3}
By incorporating the bisection-deflation and trust region algorithms above, we obtain the IAOBDM for multiple solutions of \eqref{20240110eq1.1}, whose details are presented in the following {Algorithm}.

\begin{center}
\begin{breakablealgorithm}
\setstretch{1.2}
\caption{IAOBDM for multiple solutions of \eqref{20240110eq1.1}}\label{2024080903}
\begin{algorithmic}[0]
\State {\bf Step 1.}\vspace{-1cm}
\State
\begin{adjustwidth}{1.5cm}{0.5cm}
Solve \eqref{20240112eq2.2} or \eqref{eq2.2.4.2} \textbf{Algorithm} \ref{2024080902} to obtain solution $\hat{u}^{(0)}$. Then set $\alpha_{0} = \|\hat{u}^{(0)}\|_{2}$ and $\chi_{0}= \hat{u}^{(0)} / \alpha_{0}$,
		such that $\chi_{0}$ is a unit basis function.
\end{adjustwidth}
\State {\bf Step 2.}\vspace{-1cm}
\State
\begin{adjustwidth}{1.5cm}{0.5cm}
Let $\tilde{u}^{(0)} = \alpha \chi_{0}$, and substitute it into the original equation \eqref{20240111eq2.4}. Based on the Spectral-Galerkin method, a nonlinear algebraic equation on $\alpha$ can be obtained, denoted by $\omega(\alpha)$. With the \textbf{Algorithm} \ref{2024080901}, we can obtain multiple values of $\alpha$, and denote them by $\{\alpha^{(0)}_{0}, \alpha^{(1)}_{0}, \alpha^{(2)}_{0}, \cdots\}$, which leads to
		$\tilde{u}^{(0)}_{i} = \alpha^{(i)}_{0} \chi_{0} ~ (i = 0, 1, 2, \cdots)$.
\end{adjustwidth}
\State {\bf Step 3.}\vspace{-1cm}
\State
\begin{adjustwidth}{1.5cm}{0.5cm}
To obtain multiple solutions with higher precision, we consider $\hat{u}_{1} = \alpha_{1, 0}\chi_{0} + \alpha_{1, 1}\chi_{1}$, where $\alpha_{1, 0}, \alpha_{1, 1}$ and the basis function $\chi_{1}$ are unknown variables to be solved. To determine them, the following system is solved by  \textbf{Algorithm} \ref{2024080902}.\\
		\begin{equation}\label{eq3.1}
		\begin{cases}
		\frac{{\bs F}(\hat{u}_{1})}{\alpha_{1, 1}} = {\bs 0}, \quad \alpha_{1,1} \neq 0,\\
		(\chi_{0},~ \chi_{1}) = 0,  \\
		(\chi_{1},~ \chi_{1}) = 1,
		\end{cases}
		\end{equation}
		where $(\cdot,~ \cdot)$ represents the inner product. Consequently, we can obtain multiple $\{\alpha_{1, 0}, \alpha_{1, 1}, \chi_{1}\}$ and denote them by $\bigl\{\alpha_{1, 0}^{(j)}, \alpha_{1, 1}^{(j)}, \chi_{1}^{(j)}\bigr\}~ (j = 1, 2, \ldots, k)$.
\end{adjustwidth}
\State {\bf Step 4.}\vspace{-1cm}
\State
\begin{adjustwidth}{1.5cm}{0.5cm}
Based on the Gram-Schmidt orthogonalization method, a sequence of adaptive orthogonal basis functions are generated from $ \{\chi_0, \chi_{1}^{(1)}, \chi_{1}^{(2)}, \ldots, \chi_{1}^{(k)}\}$ obtained in \textbf{Step 1} and \textbf{Step 3}, and we denote them by $S = \text{span} \big\{\chi_{0}^{(1)}, \chi_{1}^{(1)}, \ldots, \chi_{k}^{(1)} \big\}$.
\end{adjustwidth}
\State {\bf Step 5.}\vspace{-1cm}
\State
\begin{adjustwidth}{1.5cm}{0.5cm}
Let $\hat{u}^{(1)} = \hat{\alpha}_{1, 0}\chi^{(1)}_{0} + \hat{\alpha}_{1, 1}\chi^{(1)}_{1}$, and substitute it into \eqref{20240111eq2.4}. when $\hat{\alpha}_{1, 0} := \alpha_{0}^{(i)}~(i = 0, 1, 2, \cdots)$ is fixed, a nonlinear algebraic equation on $\hat{\alpha}_{1, 1}$ can be derived with the Spectral-Galerkin method. Similar to \textbf{Step 2}, multiple values of $\hat{\alpha}_{1, 1}$ can be obtained by \textbf{Algorithm} \ref{2024080901}, and denoted by $\{\hat{\alpha}^{(i, 0)}_{1, 1}, \hat{\alpha}^{(i, 1)}_{1, 1}, \cdots, \hat{\alpha}^{(i, m)}_{1, 1}\}$, where $m \geq 1$.
\end{adjustwidth}
\State {\bf Step 6.}\vspace{-1cm}
\State
\begin{adjustwidth}{1.5cm}{0.5cm}
Based on the Spectral-Galerkin method, $\tilde{u}^{(1)} := \hat{\alpha}_{1, 0}\chi^{(1)}_{0} + \hat{\alpha}_{1, 1}\chi^{(1)}_{1}$ is substituted into \eqref{20240111eq2.4}, where the unknown variables are $\hat{\alpha}_{1, 0}$ and $\hat{\alpha}_{1, 1}$. \textbf{Algorithm} \ref{2024080902} is used to obtain them with initial guesses $\{\alpha^{(i)}_{0}, \hat{\alpha}^{(i, n)}_{1, 1}\}~(n = 0, 1, \cdots, m)$. Consequently, multiple $\{\hat{\alpha}_{1, 0},~ \hat{\alpha}_{1, 1}\}$ can be obtained, and denoted by $\{\tilde{\alpha}^{(p)}_{1, 0},~ \tilde{\alpha}^{(p)}_{1, 1}\}~(p = 0, 1, \ldots)$.\vspace{0.1cm}
\end{adjustwidth}
\State {\bf Step 7.}\vspace{-1cm}
\State
\begin{adjustwidth}{1.5cm}{0.5cm}
To obtain more solutions, we consider $\tilde{u}^{(1)} = \tilde{\alpha}^{(p)}_{1, 0}\chi^{(1)}_{0} + \tilde{\alpha}^{(p)}_{1, 1}\chi^{(1)}_{1} + \hat{\alpha}_{1, 2}\chi^{(1)}_{2}$, where $\hat{\alpha}_{1, 2}$ is a unknown variable to be solved. Similarly to the process from \textbf{Step 5} to \textbf{Step 6}, more solutions can be obtained until the adaptive basis functions in \textbf{Step 4} are fully used.
\end{adjustwidth}
\end{algorithmic}
\end{breakablealgorithm}
\end{center}
\vspace{0.15cm}
{\color{black}Some remarks are listed in order.}

\begin{rem}\label{rem1.1}
The functions used to construct adaptive bases are solved by 
\textbf{Algorithm} \ref{2024080902} in \textbf{Step 1} and \textbf{Step 3},
which are the most expensive parts of the overall algorithm. The trust region method (\textbf{Algorithm} \ref{2024080902}) is adopted since it allows flexible initial values and, at the same time, has {\color{black}good} convergence.
\end{rem}

\begin{rem}\label{rem3.3} 
   {\color{black} In {\textbf{Step 2} and \textbf{Step 5}}, scalar NLAS are solved by using
    \textbf{Algorithm} \ref{2024080901}. It is worth pointing out that in the previous study \cite{Li2024}, the adaptive orthogonal basis method is only suitable for the polynomial nonlinearities, which greatly limits its scope of application. In this paper, the bisection-deflation algorithm presented in \textbf{Algorithm} \ref{2024080901} is introduced to provide good initial guesses for finding multiple solutions to more general nonlinear differential equations}.

\end{rem}

\begin{rem}\label{rem3.4} 
    In {\textbf{Step 6} and \textbf{Step 7}}, {\color{black}when more than} one basis functions are used in the refinement procedure, the resulting NLASs have more than one equation, {\color{black}thus we use  \textbf{Algorithm} \ref{2024080902}, instead of 
    \textbf{Algorithm} \ref{2024080901}. Since the dimensions of these NLASs are small, their computational costs are much less than solving \eqref{20240112eq2.2} or \eqref{eq2.2.4.2} for multiple solutions directly using the deflation method.} 
\end{rem}





\section{Numerical results}\label{sect4}
\begin{figure}[!h]
\begin{center}	\includegraphics[width=5.5cm,height=4.6cm]{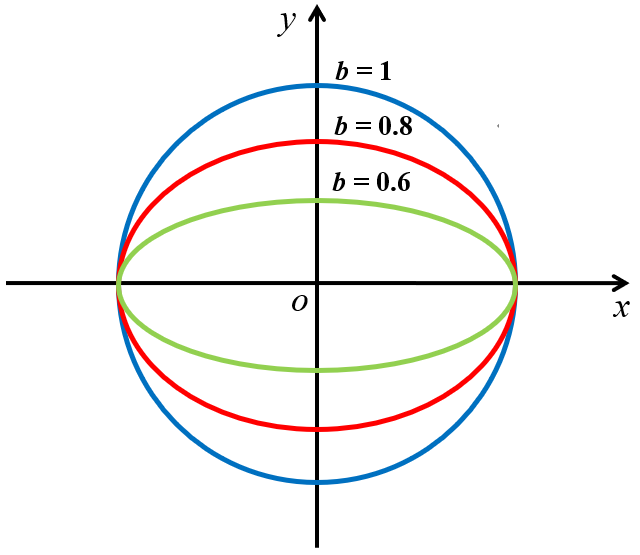}
\caption{Different regions with varying $b$.}\label{figure1}
\end{center}
\end{figure}

\setcounter{thm}{0}
In this section, we first validate the accuracy and efficiency of the proposed IAOBDM described in \textbf{Algorithm} \ref{2024080903}, then apply it to some practical problems. In subsequent examples, we consider multiple solutions on the unit disk $\Omega = \{(x, y)| x^2 + y^2 \leq 1\}$ and the elliptic region $\hat{\Omega} = \{ (x, y)| x^2 + \frac{y^2}{b^2} \leq 1, 0< b < 1\}$ (see Fig.\ref{figure1}). The corresponding variational functional values of these multiple solutions will be computed with varying $b$. Some theoretical results on multiple solutions of these examples will be numerically verified. All calculations are carried out on a workstation with Intel(R) Core(TM) i9-13900HX (2.20GHz) and 120GB RAM.

{\em \textbf{Example 1.}} Consider the following elliptic sine-Gordon equation~\cite{Chen2003ANO}:
\begin{equation}\label{2024081801}
\triangle u + \lambda \sin u = 0,   \quad\quad   x = (x_1, x_2)\in \Omega, \quad\quad  \lambda > 0,
\end{equation}
with homogeneous Dirichlet boundary condition 
\begin{equation}\label{eq:DirichletBC}
    u({\bs x}) = 0, \quad \text{on }\partial \Omega,
\end{equation}
or homogeneous Neumann boundary condition
\begin{equation}\label{eq:NeumannBC}
    \partial_{\bs n} u({\bs x}) = 0, \quad \text{on }\partial \Omega.
\end{equation}

A variational {\color{black}functional} corresponding to \eqref{2024081801}-\eqref{eq:NeumannBC} is
\begin{equation}\label{2024081802}
J(u) = \int_{\Omega}\bigl[\frac{1}{2}|\nabla u|^2 - \lambda(1 - \cos u)\bigr]dx,  \quad\quad\quad  u \in H^1(\Omega),
\end{equation}
where $H^1(\Omega)$ denotes the usual Sobolev space of order 1. If Dirichlet boundary condition \eqref{eq:DirichletBC} is used, we replace $H^1(\Omega)$ by $H^1_0(\Omega)$.

\hspace{-0.6cm}\underline{\textbf{Case 1.}} (\textbf{Dirichlet boundary condition}) 
Let $\{ \hat{\lambda}_j \mid 0 = \hat{\lambda}_1 < \hat{\lambda}_2 < \cdots, j = 1, 2, \ldots \}$ be the set of eigenvalues of the following problem
\begin{equation}
\begin{cases}
-\triangle\phi_j = \lambda_j\phi_j  \quad\quad\quad  \textrm{in}\;\Omega,\\
\phi_j = 0   \quad\quad\quad \quad\quad\;\;\,\;\;\textrm{on}\;\partial\Omega.
\end{cases}
\end{equation}
The following theoretical result on multiple solutions of \eqref{2024081801}-\eqref{eq:DirichletBC} has been shown by~\cite{Chen2003ANO}:
\begin{thm}
\eqref{2024081801}-\eqref{eq:DirichletBC} has at least four nontrivial solutions for $\lambda > \hat{\lambda}_2$.
\end{thm}

\begin{figure}[h]
\begin{center} \subfigure[I]{\includegraphics[width=3.35cm,height=3.3cm]{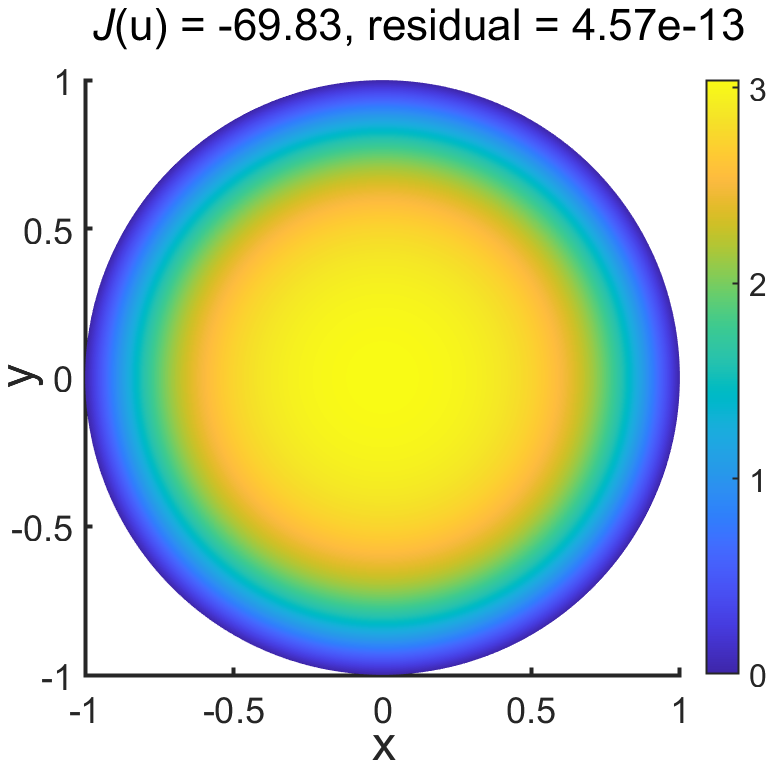}}\quad \subfigure[II]{\includegraphics[width=3.35cm,height=3.3cm]{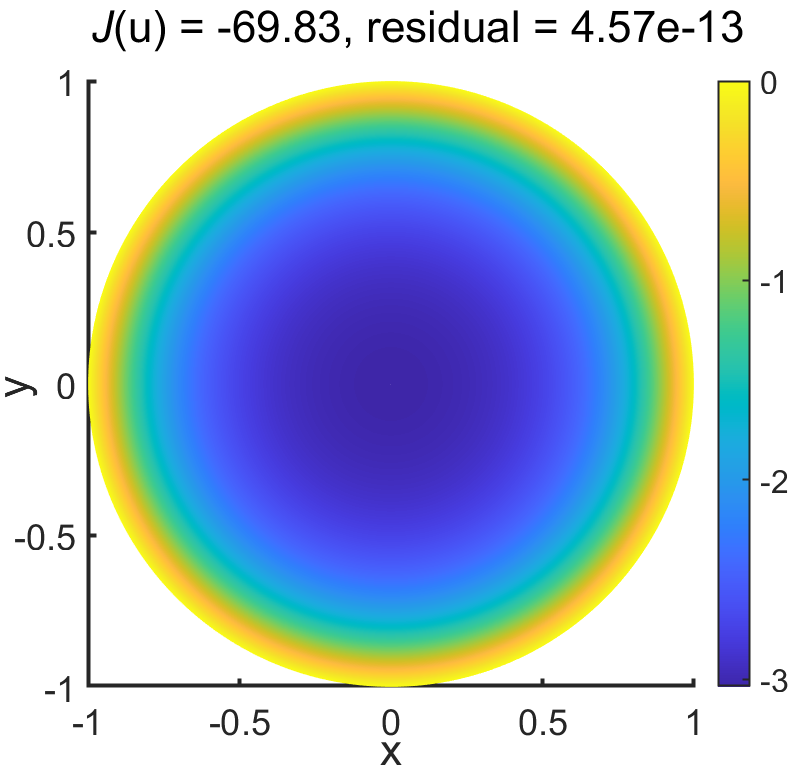}}\quad \subfigure[III]{\includegraphics[width=3.35cm,height=3.3cm]{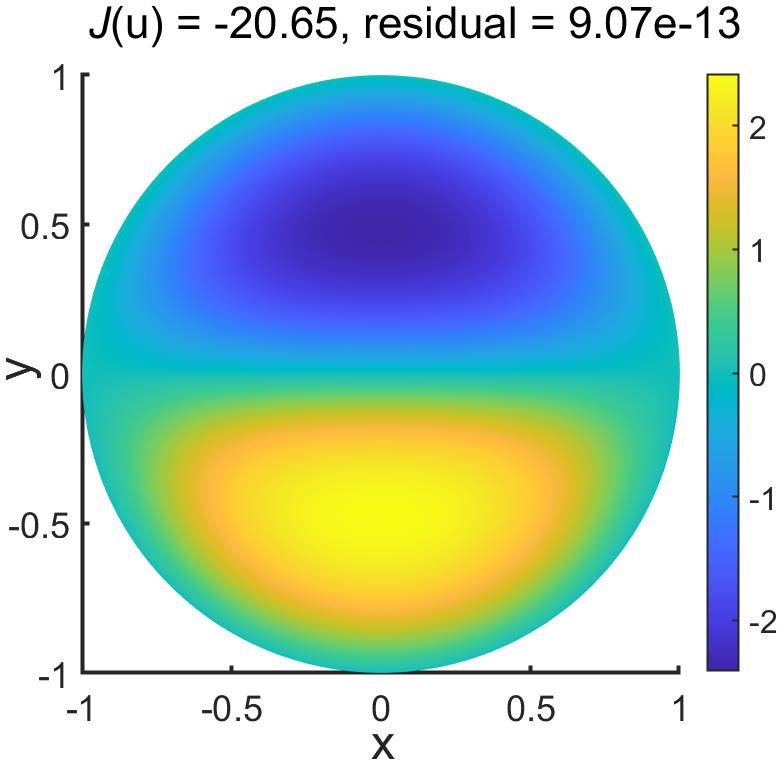}}\quad \subfigure[IV]{\includegraphics[width=3.35cm,height=3.3cm]{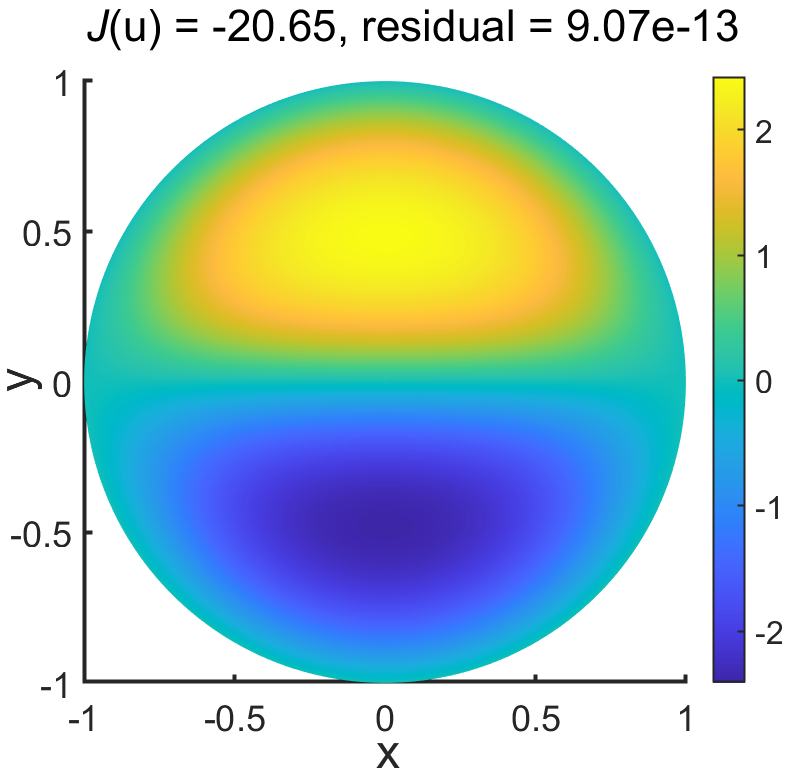}}
\caption{Multiple solutions of \eqref{2024081801}-\eqref{eq:DirichletBC} with $\lambda=30$.}\label{example01}
\end{center}
\end{figure}

\begin{figure}[H]
\begin{center}	 \subfigure[$b = 0.9$]{\includegraphics[width=2.9cm,height=3.3cm]{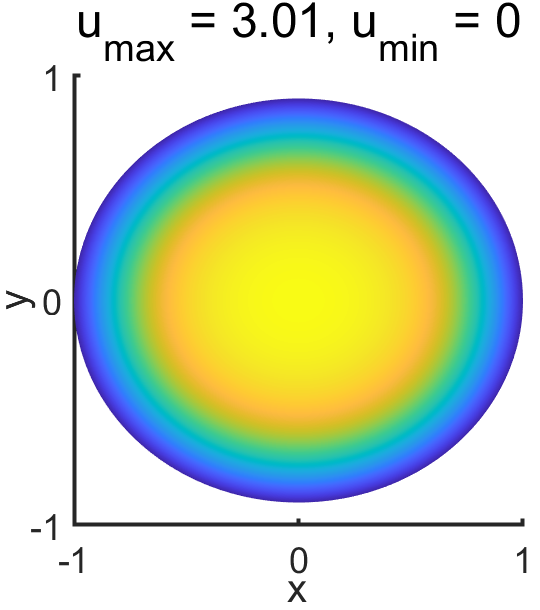}}\quad\quad \subfigure[$b = 0.8$]{\includegraphics[width=3cm,height=3.3cm]{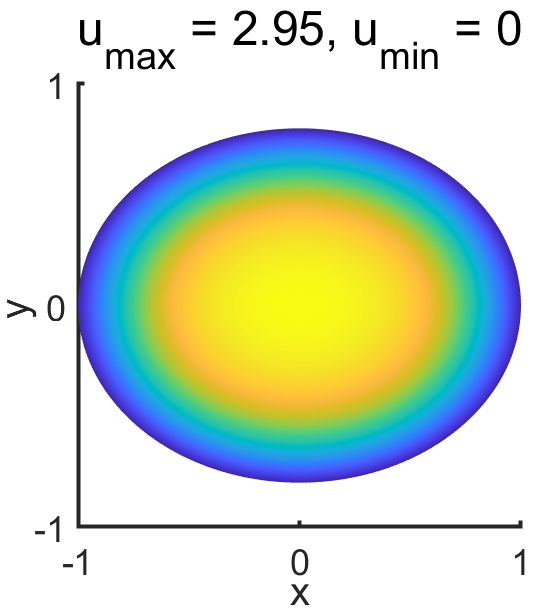}}\quad\quad \subfigure[$b = 0.7$]{\includegraphics[width=3cm,height=3.3cm]{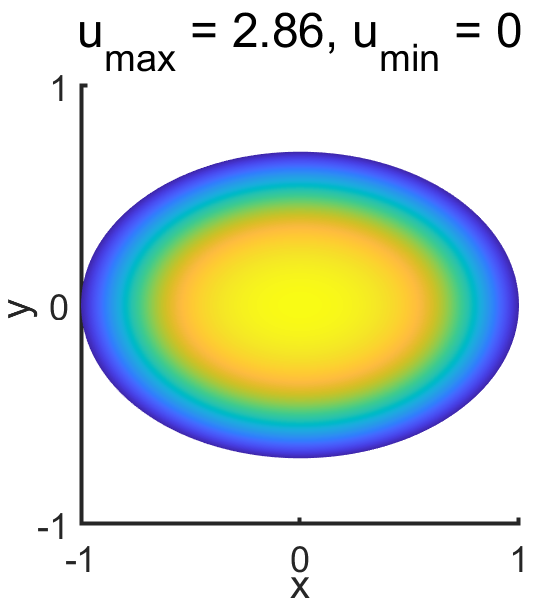}}\quad\quad 
\caption{The type-I solution to \eqref{2024081801}-\eqref{eq:DirichletBC} with varying $b$.}\label{2024101401example23}
 \subfigure[$b = 0.9$]{\includegraphics[width=3cm,height=3.3cm]{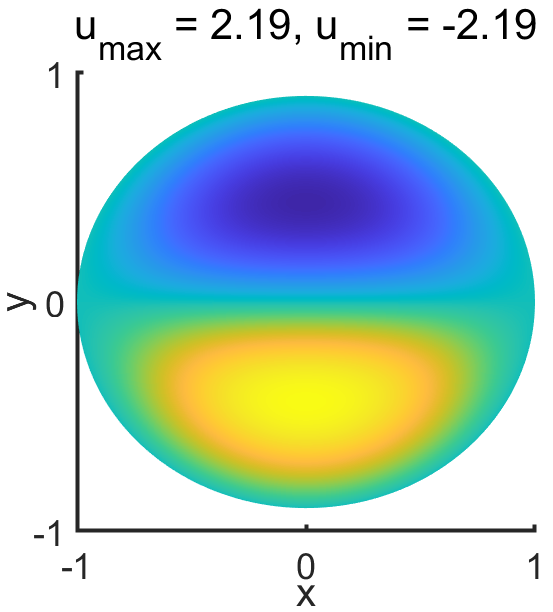}}\quad\quad \subfigure[$b = 0.8$]{\includegraphics[width=3cm,height=3.3cm]{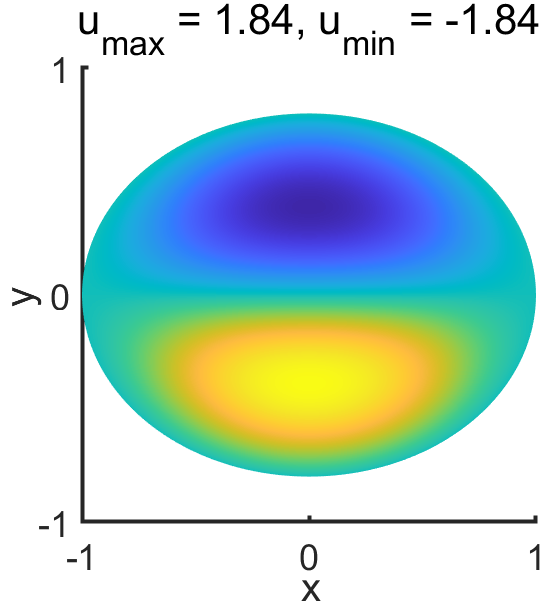}}\quad\quad \subfigure[$b = 0.7$]{\includegraphics[width=2.9cm,height=3.25cm]{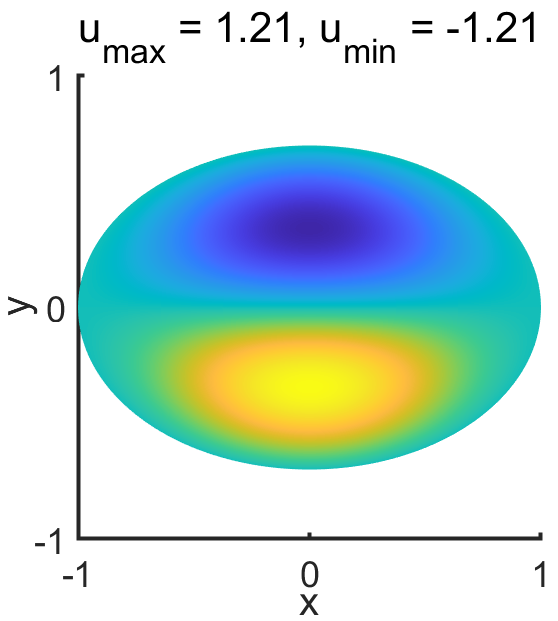}}\quad\quad
\caption{The type-III solution to \eqref{2024081801}-\eqref{eq:DirichletBC} with varying $b$.}\label{2024101402example23}
\end{center}
\end{figure}

By the Legendre--Galerkin method, we find that $\hat{\lambda}_2 \approx 14.22$. As shown in Fig.\ref{example01}, four nontrivial solutions are obtained, where the type-I and type-II solutions are rotationally symmetric about the axis $r = 0$. Although type-III and type-IV solutions shown in Figs.\ref{example01}(c)-\ref{example01}(d) are not rotationally invariant, they are still solutions after an arbitrary rotation about $r=0$. The residuals of these multiple solutions are marked in the title text of the figures, verifying that these solutions are of high accuracy.
Moreover, the type-I and type-III solutions of \eqref{2024081801} with varying $b$ are presented in Figs.\ref{2024101401example23}-\ref{2024101402example23}. In addition, $J(u)$ defined in \eqref{2024081802} with $\lambda =30$ are plotted with varying $b$ in Fig.\ref{example04}(a). 

\hspace{-0.6cm}\underline{\textbf{Case 2.}} (\textbf{Neumann boundary condition}) Let $\{ \lambda_j \mid 0 = \lambda_1 < \lambda_2 < \cdots, j = 1, 2, \ldots \}$ be the set of eigenvalues of the following problem
\begin{equation}
\begin{cases}
-\triangle\psi_j = \lambda_j\psi_j  \quad\quad\quad  \textrm{in}\;\Omega,\\
\frac{\partial\psi_j}{\partial n} = 0   \quad\quad\quad \quad\quad\;\;\, \textrm{on}\;\partial\Omega.
\end{cases}
\end{equation}
Some theoretical results have been given as follows~\cite{Chen2003ANO}:
\begin{thm}\label{2024082901}
If $\lambda > \lambda_2$, \eqref{2024081802}-\eqref{eq:NeumannBC} has at least two nonconstant critical points $u$ in $H^1(\Omega)$. Moreover, $u+2k\pi$ are also solutions of \eqref{2024081801}-\eqref{eq:NeumannBC} for any integer $k$.
\end{thm}
\begin{thm}\label{2024082902}
For any $\lambda > 0$ and integer $k$, \eqref{2024081801}-\eqref{eq:NeumannBC} does not admit any nonconstant solution $u$ such that $0 < u < \pi$, and it does not admit any nonconstant solution $u$ such that $2k\pi < u < (2k+1)\pi$.\vspace{-0.25cm}
\end{thm}

\begin{figure}[h]
\begin{center} \subfigure[I]{\includegraphics[width=4.15cm,height=3.9cm]{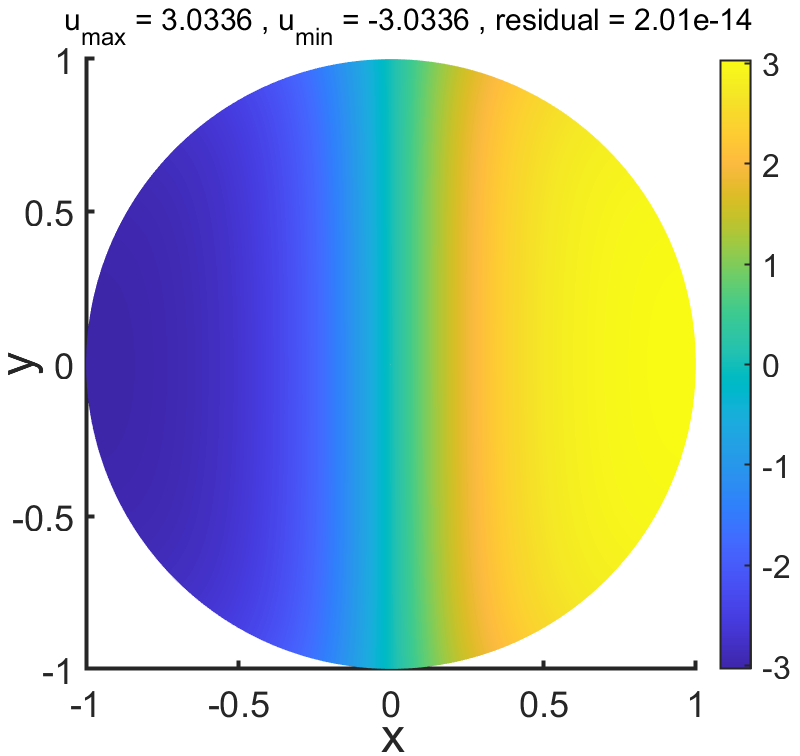}}\quad\quad \subfigure[I+2$\pi$]{\includegraphics[width=4cm,height=3.9cm]{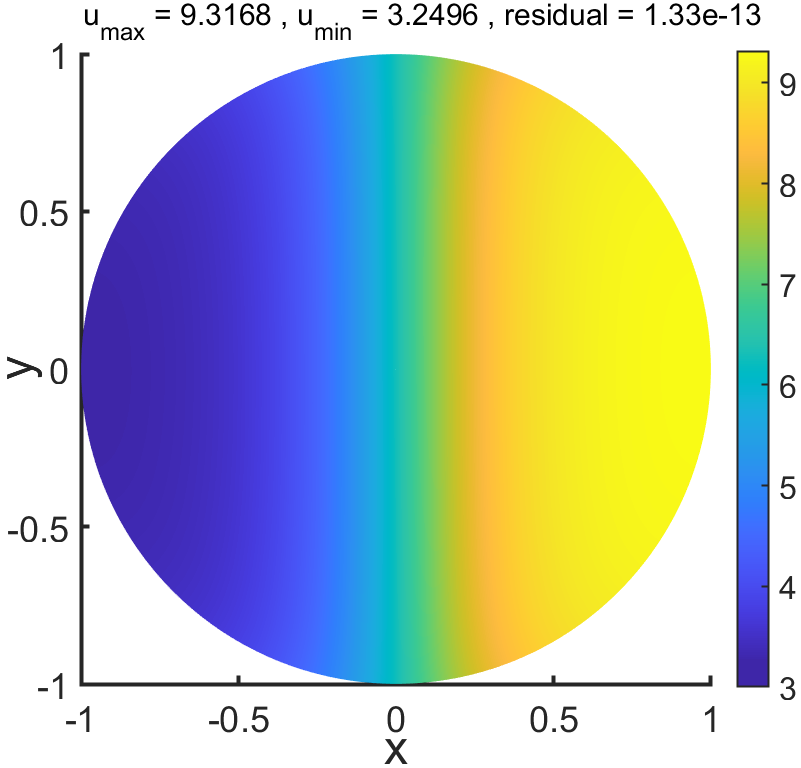}}\quad\quad
\subfigure[I+4$\pi$]{\includegraphics[width=4.1cm,height=3.9cm]{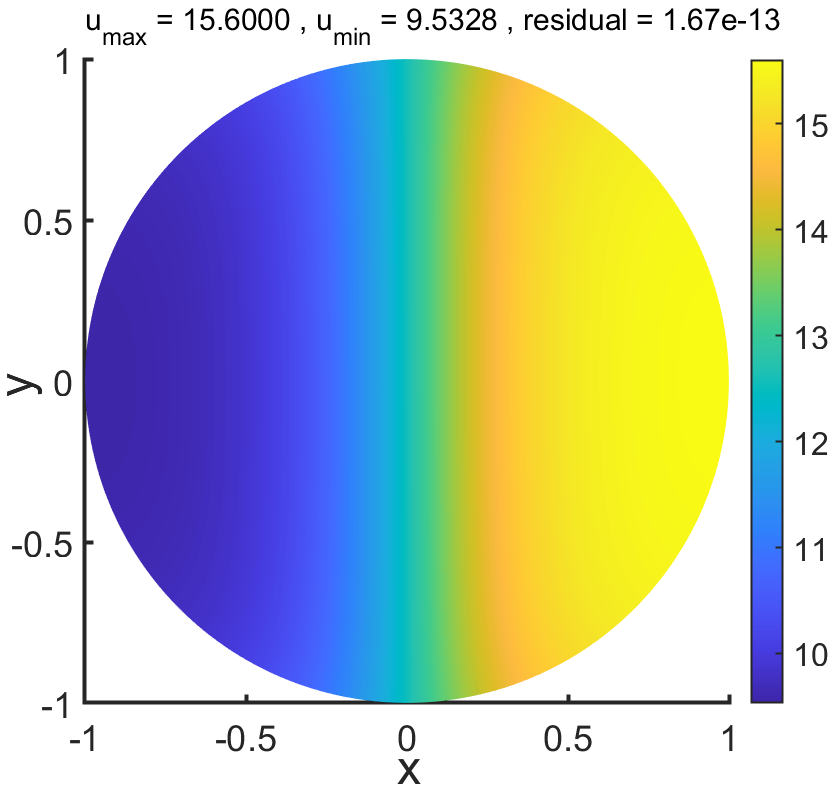}}\\
\subfigure[II]{\includegraphics[width=4.1cm,height=3.9cm]{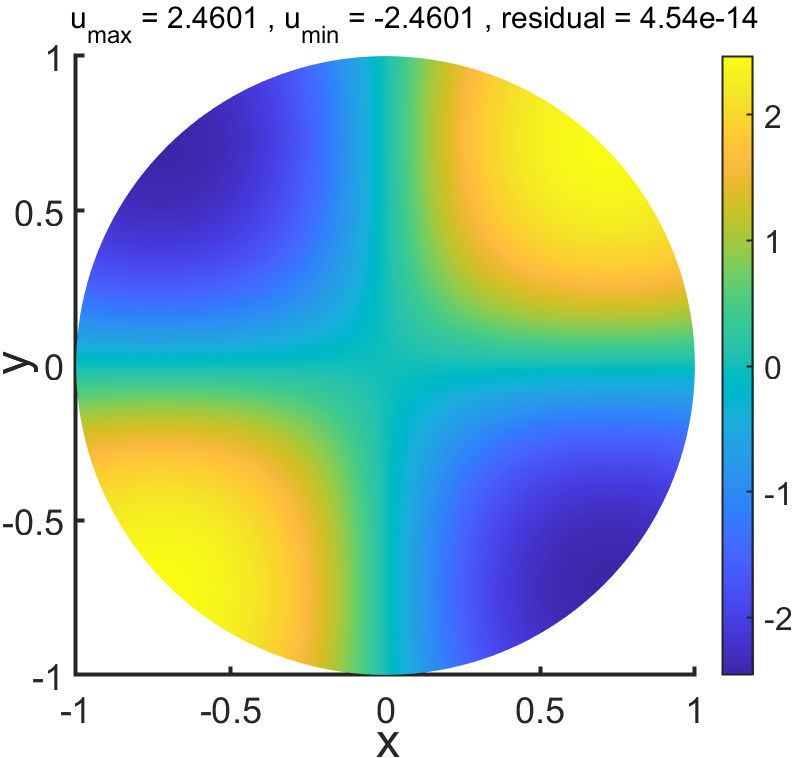}}\quad\quad \subfigure[II+2$\pi$]{\includegraphics[width=4.1cm,height=3.9cm]{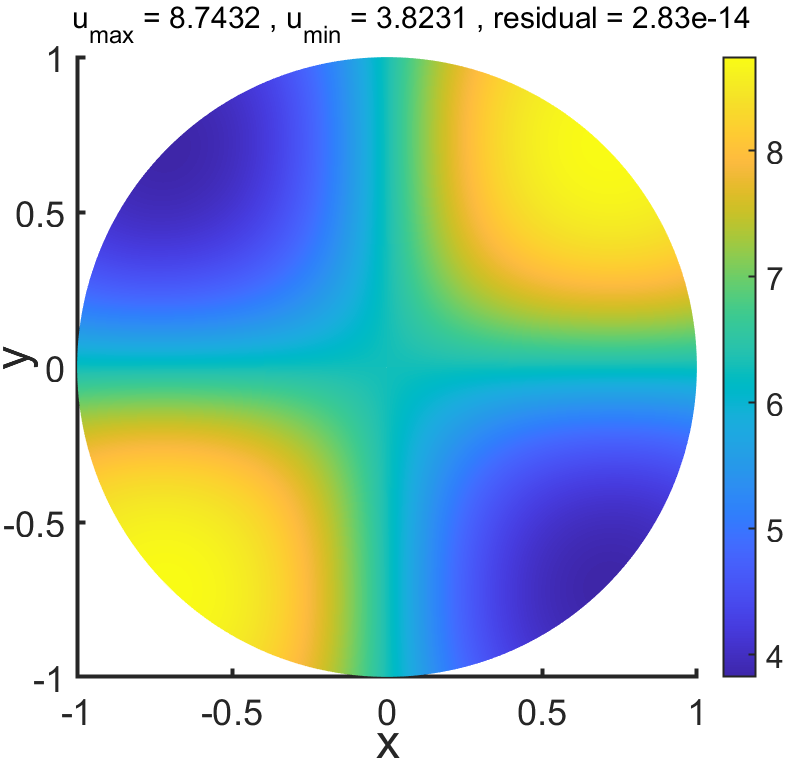}}\quad\quad
\subfigure[II+4$\pi$]{\includegraphics[width=4.1cm,height=3.9cm]{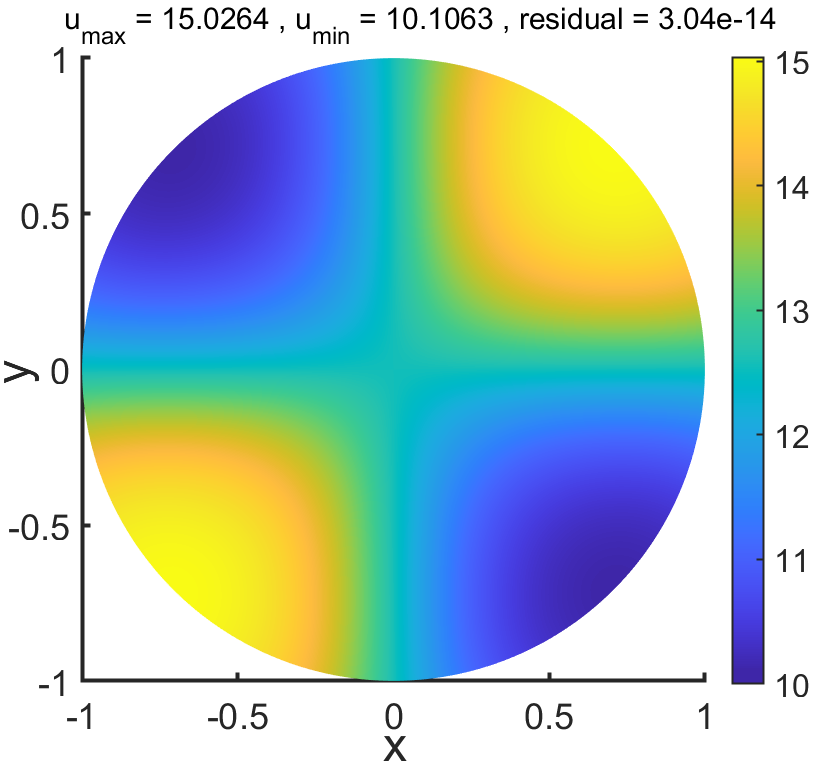}}
\caption{Multiple solutions of \eqref{2024081801}-\eqref{eq:NeumannBC} with $\lambda = 20$.}\label{example03}\vspace{-0.7cm}
\end{center}
\end{figure}
By the Legendre--Galerkin method, $\lambda_2 \approx 9.43$ is computed. As shown in Fig.\ref{example03}, two nonconstant critical points $u$ (i.e., Figs.\ref{example03}(a)-\ref{example03}(d)) are obtained, which verifies the result presented in \textbf{Theorem} \ref{2024082901}. Moreover, $u+2k\pi$ can be easily verified as the solutions of \eqref{2024081801}-\eqref{eq:NeumannBC} for \textbf{case 2}. Meanwhile, the solutions presented in Figs.\ref{example03}(b) and \ref{example03}(c) and Figs.\ref{example03}(e) and \ref{example03}(f) agree well with the result in \textbf{Theorem} \ref{2024082902}. In addition, $J(u)$ with $\lambda = 20$ are plotted in Fig.\ref{example04}(b) with varying $b$ for the type-I and type-II solutions (see Figs.\ref{example1}-\ref{example2}) to \eqref{2024081801}-\eqref{eq:NeumannBC}.

{\color{black}In addition, a comparison of different methods is presented in Table \ref{NNewTable61}. The residual $\|{\bs F}\|_{\infty}$ is shown with some increasing numbers of iterations (i.e., $n_{it}$) for homogeneous Dirichlet and Neumann boundary conditions (i.e., \textbf{Case 1} and \textbf{Case 2}), where $\hat{\lambda}_1, \hat{\lambda}_2, \alpha^{(0)}_0$ and $\alpha^{(1)}_0$ represent different initial guesses. Please refer to~\cite{chen2004search,Li2024} and \textbf{Algorithm} 3.3 for more details. To use} 
\begin{figure}[!h]
\begin{center} 
\subfigure[Dirichlet boundary with $\lambda = 30$]{\includegraphics[width=4.9cm,height=4.3cm]{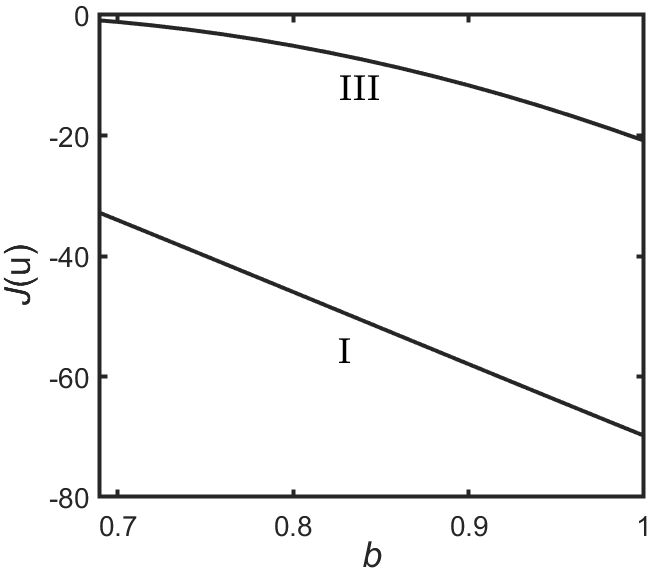}} \quad\quad\quad
\subfigure[Neumann boundary with $\lambda = 20$]{\includegraphics[width=4.9cm,height=4.3cm]{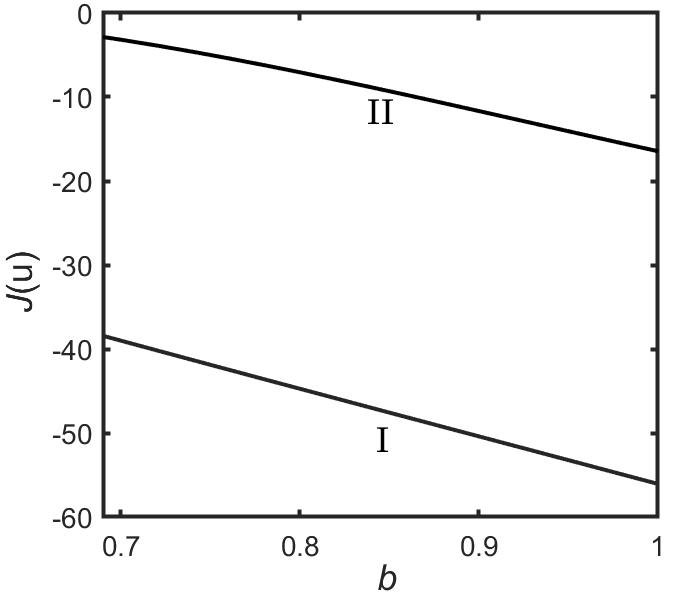}}
\caption{The dependence of $J(u)$ defined in \eqref{2024081802} on geometry parameter $b$ for different solutions.}\label{example04}
\end{center}\vspace{-0.3cm}
\end{figure}

\begin{figure}[H]
\begin{center}	 \subfigure[$b = 0.9$]{\includegraphics[width=3cm,height=3.3cm]{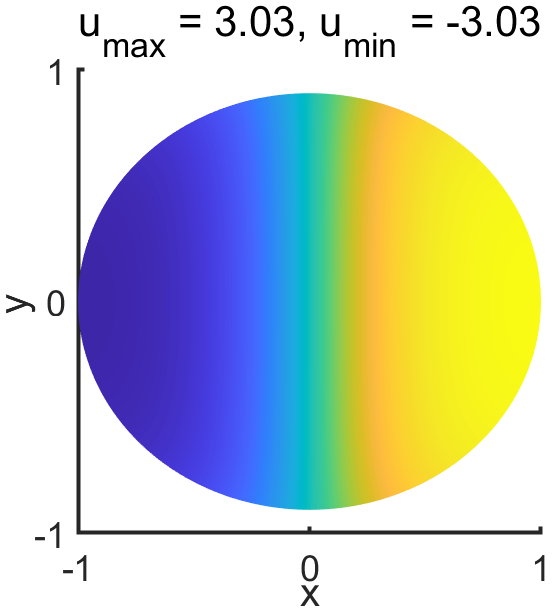}}\quad\quad \subfigure[$b = 0.8$]{\includegraphics[width=3cm,height=3.3cm]{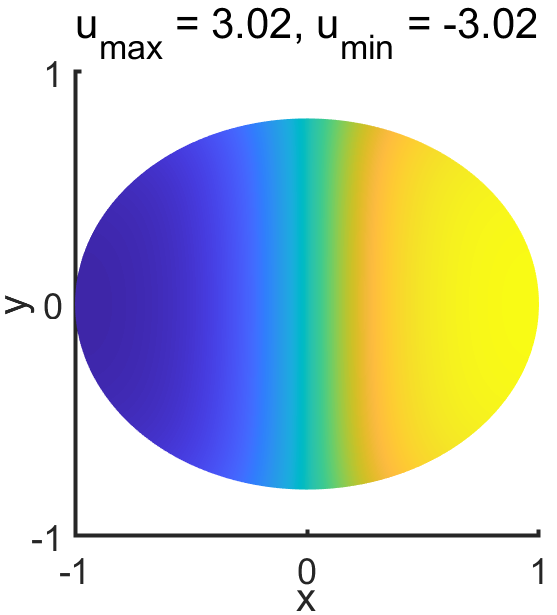}}\quad\quad \subfigure[$b = 0.7$]{\includegraphics[width=3cm,height=3.3cm]{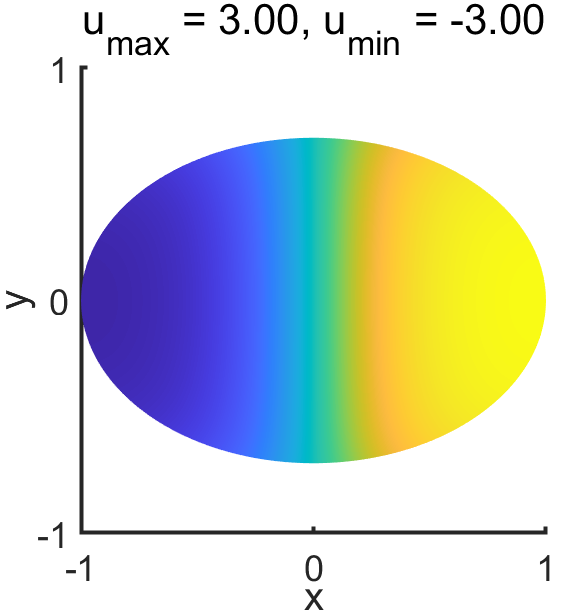}}\vspace{-0.06cm}\quad\quad 
\caption{The type-I solution to \eqref{2024081801}-\eqref{eq:NeumannBC} with varying $b$.}\label{example1}
 \subfigure[$b = 0.9$]{\includegraphics[width=3cm,height=3.3cm]{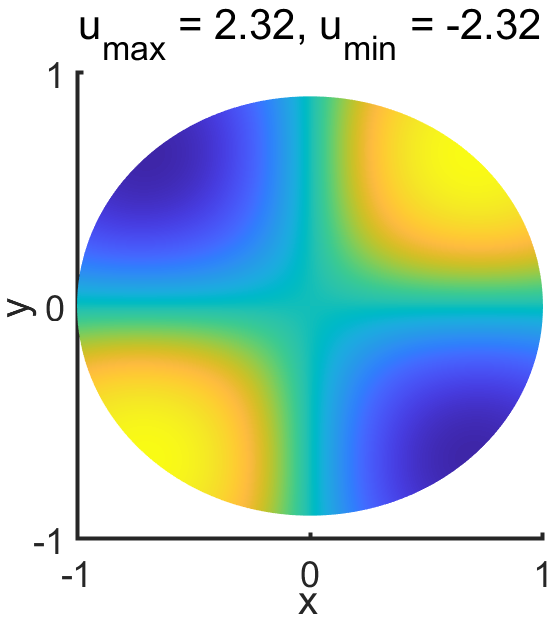}}\quad\quad \subfigure[$b = 0.8$]{\includegraphics[width=3cm,height=3.3cm]{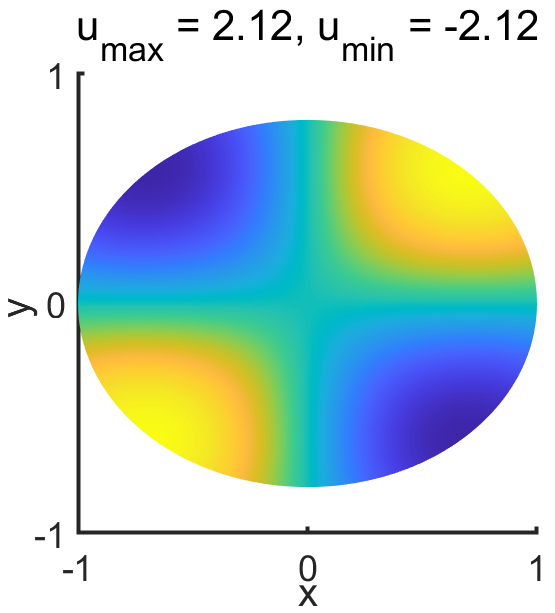}}\quad\quad \subfigure[$b = 0.7$]{\includegraphics[width=3cm,height=3.3cm]{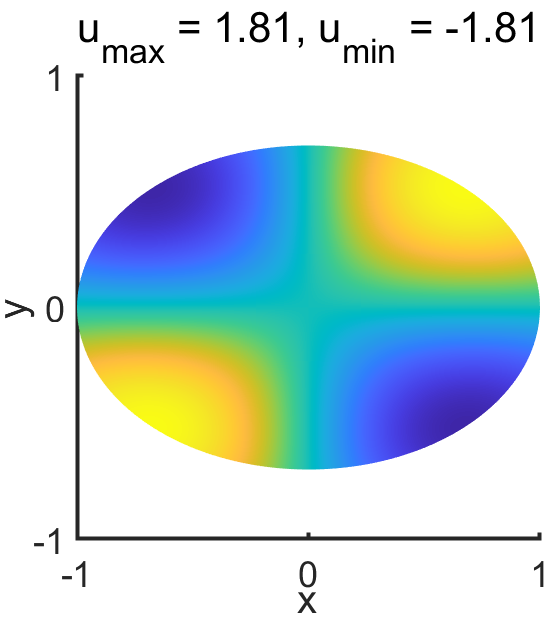}}\quad\quad
\caption{The type-II solution to \eqref{2024081801}-\eqref{eq:NeumannBC} with varying $b$.}\label{example2}\vspace{-0.9cm}
\end{center}
\end{figure}

\begin{table}[h]
\centering\small
\caption{\small A comparison of different methods. The real numbers in this table represent the residuals $\|{\bs F}\|_{\infty}$ after some given number of iterations.}
\label{NNewTable61}
\adjustbox{width=0.85\textwidth}{
\begin{tabular}{ccccccccccc}
	\hline
~ &~ &\multicolumn{3}{c}{SEM in \cite{chen2004search}}
&\multicolumn{3}{c}{method in \cite{Li2024}}
&\multicolumn{3}{c}{IAOBDM~(\textbf{Alg.} \ref{2024080903})}\\\cline{3-5}\cline{6-8}\cline{9-11}
\specialrule{0em}{1pt}{1.5pt}
~   &$\lambda$ &$n_{it}$ &$\hat{\lambda}_1$ &$\hat{\lambda}_2$
&$n_{it}$ &$\alpha^{(0)}_0$ &$\alpha^{(1)}_0$
&$n_{it}$ &$\alpha^{(0)}_0$ &$\alpha^{(1)}_0$
\\\hline   
\specialrule{0em}{1pt}{1.5pt}
\multirow{3}{*}{\textbf{Case 1}} &\multirow{3}{*}{30}  &10  &8.14e5 &1.34e4  &10  &9.70e5 &3.71e3
&10  &5.72e1 &4.85e1 \\
~  &~ &20  &9.05e3 &6.32e2  &20  &8.00e3 &1.41e2
&20  &7.92e-3 &9.59e-2\\
~ &~  &30  &1.27e1 &9.75e0  &30  &4.21e0 &9.15e-1
&30  &6.55e-7 &3.57e-5\\\hline
\specialrule{0em}{1pt}{1.5pt}
\multirow{3}{*}{\textbf{Case 2}} &\multirow{3}{*}{20}  &10  &3.76e4 &5.75e6  &10  &8.49e4 &9.34e7
&10  &6.55e1 &1.71e0 \\
~  &~ &20  &2.78e3 &6.49e4  &20  &6.78e2 &7.57e3
&20  &7.06e-4 &3.18e-5\\
~ &~  &30  &5.46e1 &1.57e1  &30  &7.43e-1 &3.92e0
&30  &2.76e-8 &4.62e-9\\\hline
\end{tabular}
}\vspace{-0.35cm}
\end{table}
\hspace{-0.55cm}{\color{black}the method in [35], we take $u-{u^3}/{6}$ to replace $\sin u$ in \eqref{2024081801}, since this method works only on polynomial nonlinearity. We observe from this table that the proposed \textbf{Algorithm} 3.3 is more efficient than the methods proposed in~\cite{chen2004search} and \cite{Li2024}.}

\clearpage
{\em \textbf{Example 2.}} Consider the following Henon equation with Dirichlet boundary\vspace{-0.2cm}
\begin{equation}\label{eq.4.1}
\begin{cases}
\triangle u + u^3=0\;\;\;\quad\quad   \textrm{in}\; \Omega,\\
u = 0   \quad\quad\quad\quad\quad\; \;\;\, \textrm{on}\; \partial\Omega.\vspace{-0.2cm}
\end{cases}
\end{equation}
The variational functional corresponding to \eqref{eq.4.1} is
\begin{equation}\label{eq22}
J(u) = \int_{\Omega}(\frac{1}{2}|\nabla u|^2 -  \frac{1}{4}u^4)dx,  \quad\quad\quad  u \in H^1_0(\Omega).\vspace{-0.2cm}
\end{equation}
\begin{figure}[h]
\begin{center}	
\subfigure[I]{\includegraphics[width=4.4cm,height=4.1cm]{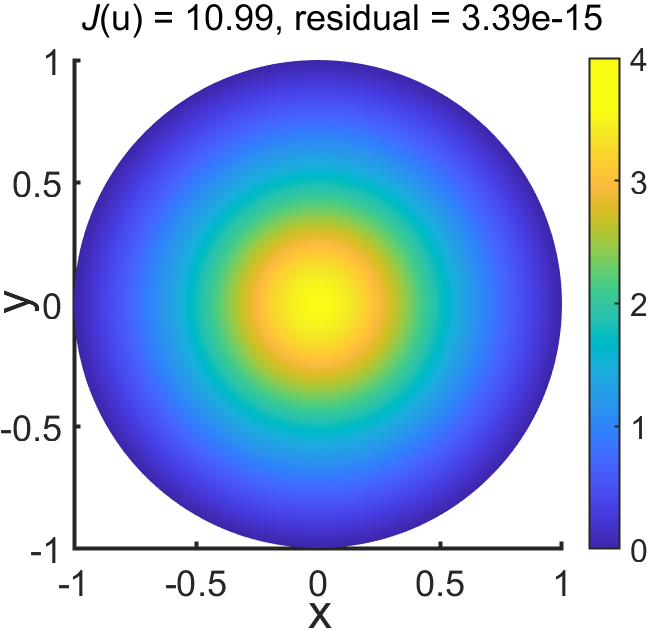}} \subfigure[II]{\includegraphics[width=4.4cm,height=4.1cm]{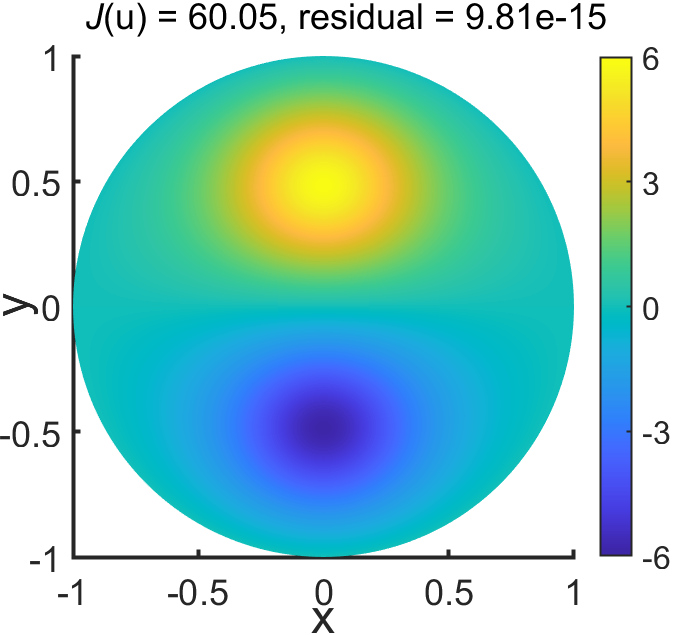}} \subfigure[III]{\includegraphics[width=4.4cm,height=4.1cm]{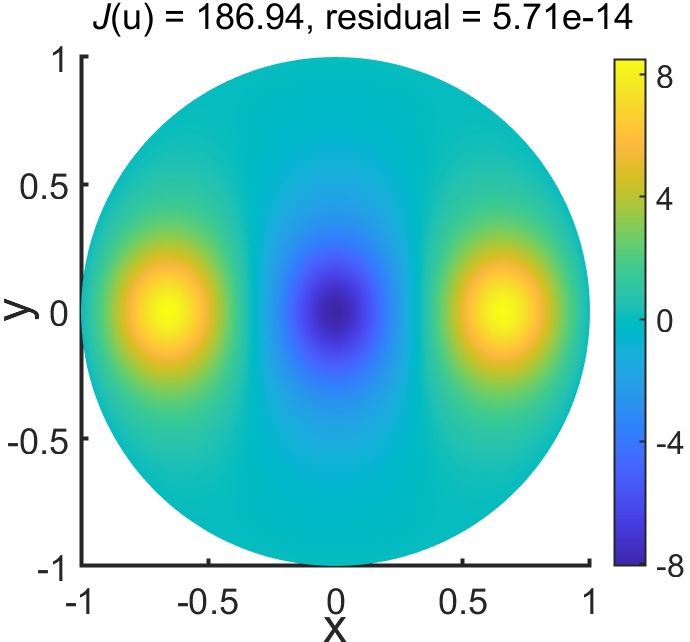}}\\ \subfigure[IV]{\includegraphics[width=4.35cm,height=4.1cm]{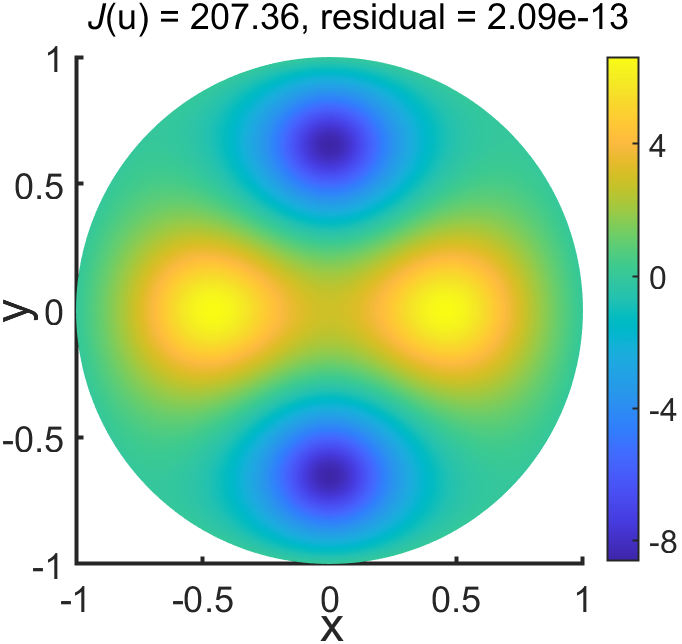}} \subfigure[V]{\includegraphics[width=4.35cm,height=4.1cm]{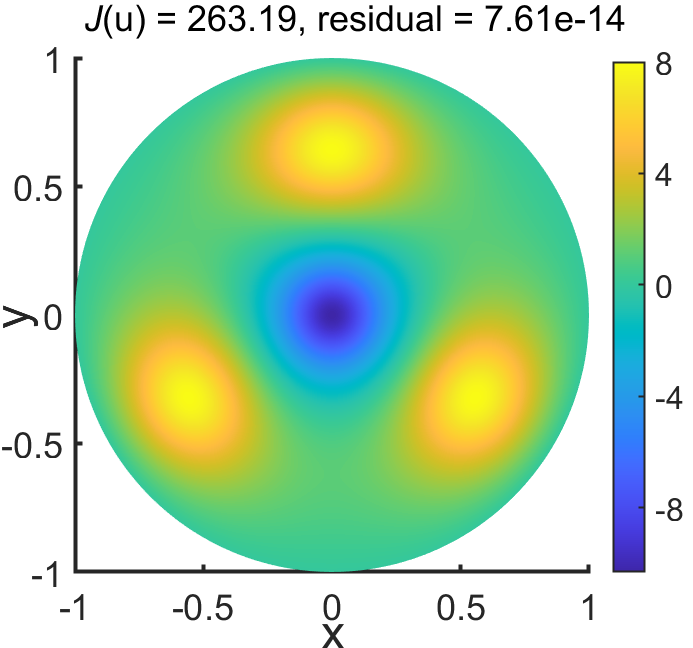}} \subfigure[VI]{\includegraphics[width=4.5cm,height=4.1cm]{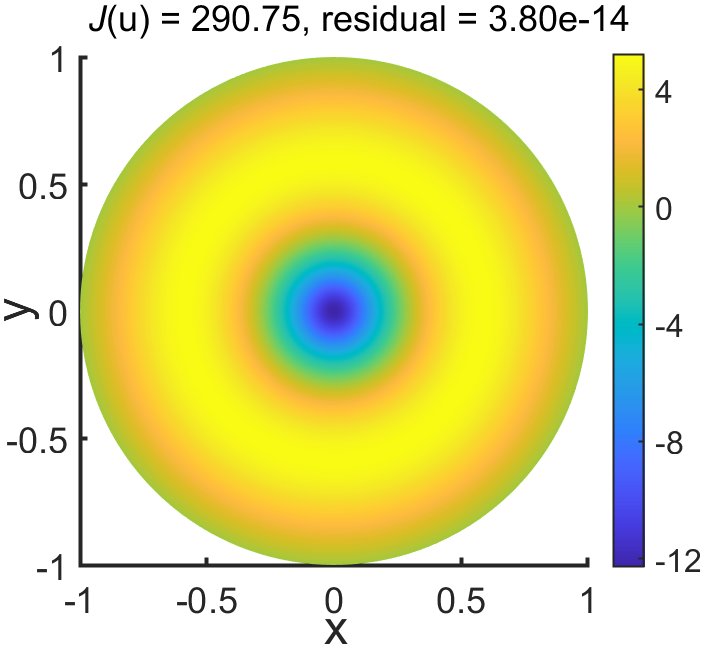}}\\ \subfigure[VII]{\includegraphics[width=4.45cm,height=4.1cm]{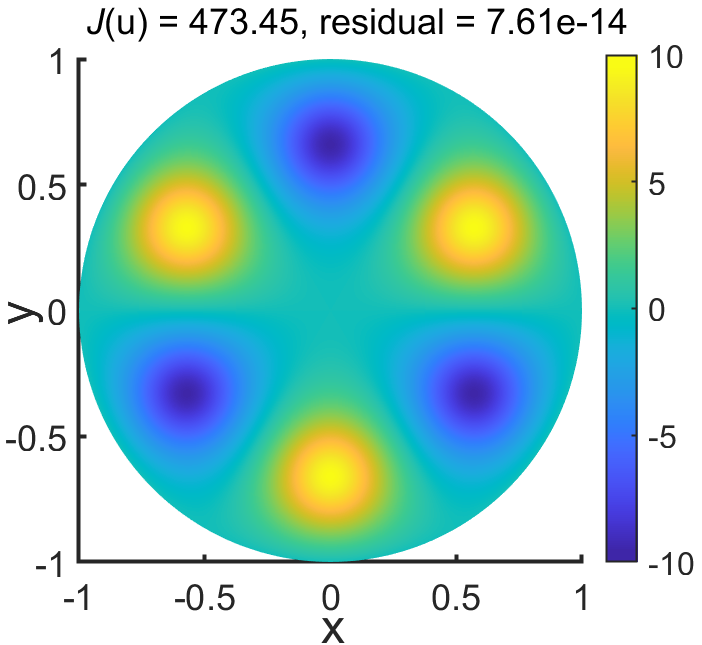}} \subfigure[VIII]{\includegraphics[width=4.45cm,height=4.1cm]{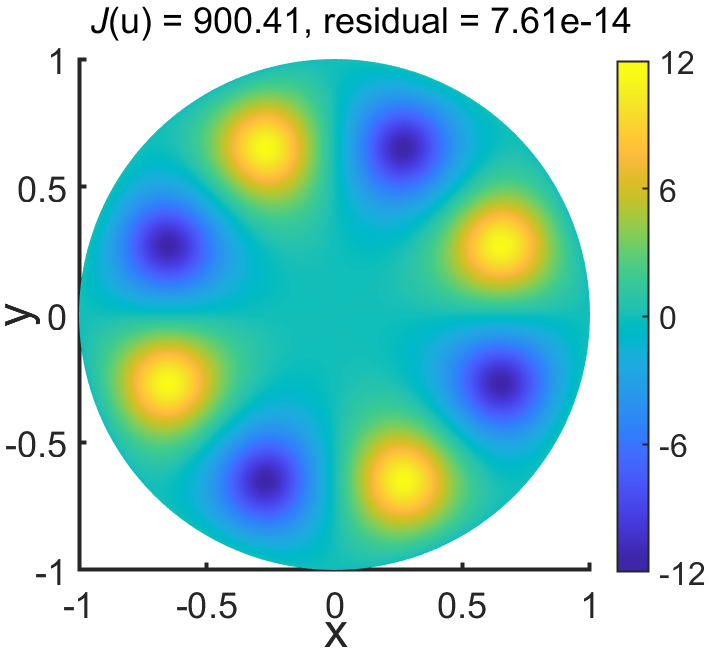}} \subfigure[IX]{\includegraphics[width=4.45cm,height=4.1cm]{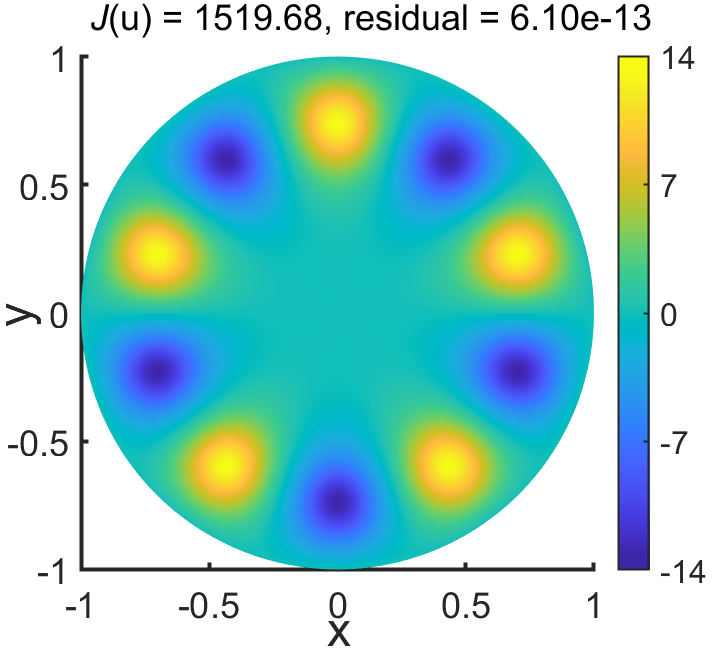}}
\caption {Multiple solutions of \eqref{eq.4.1}.}\label{example21type}
\vspace{-1cm}
\end{center}
\end{figure}

\begin{figure}[H]
\begin{center}	\subfigure{\includegraphics[width=5cm,height=5cm]{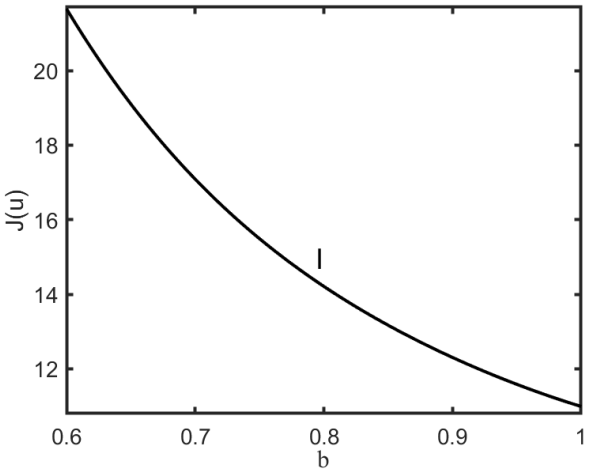}}\quad\quad\quad \subfigure{\includegraphics[width=5cm,height=5cm]{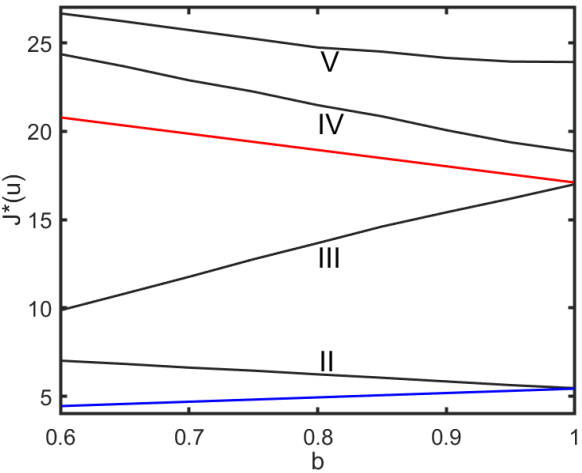}}
\caption{$J(u)$ and $J^\star(u):=J(u)/J(u_I)$ versus $b$ for multiple solutions of \eqref{eq.4.1}.}\label{20241014example22}
\end{center}
\end{figure}

\begin{figure}[H]
\begin{center}	
\subfigure[I]{\includegraphics[width=7.1cm,height=3cm]{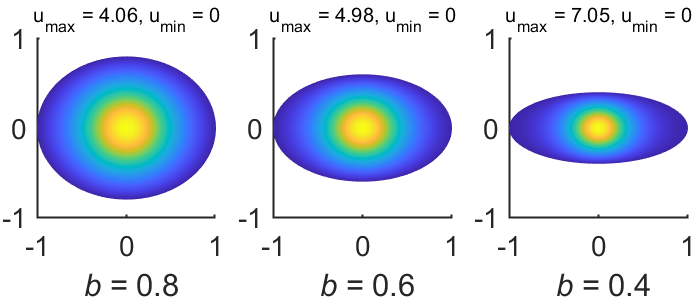}} \subfigure[II]{\includegraphics[width=7.1cm,height=3cm]{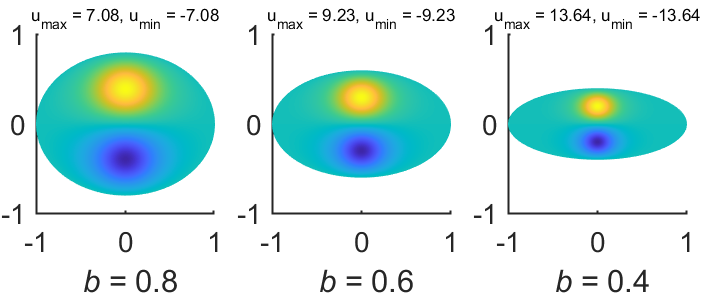}}\\ 
\subfigure[III]{\includegraphics[width=7.1cm,height=3cm]{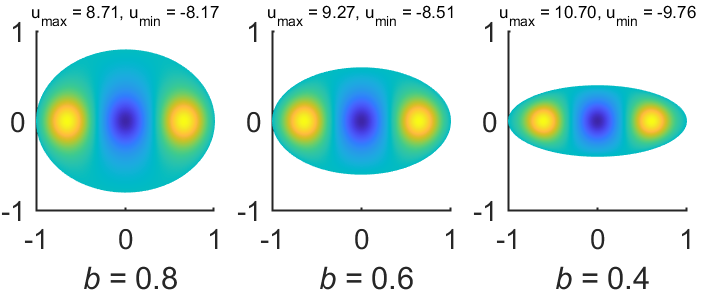}}
\subfigure[IV]{\includegraphics[width=7.1cm,height=3cm]{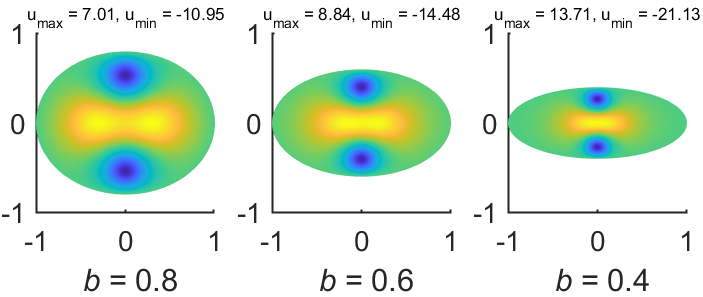}}\\ 
\subfigure[V]{\includegraphics[width=7.15cm,height=3.0cm]{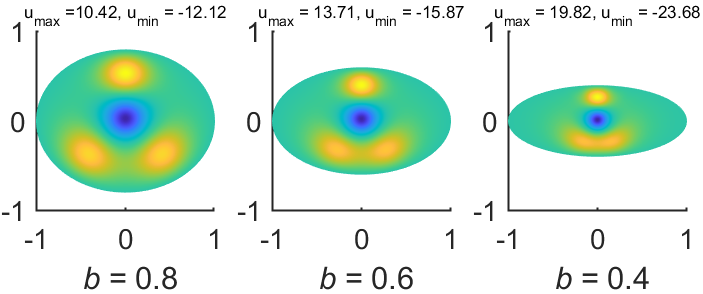}}		\subfigure[VI]{\includegraphics[width=7.15cm,height=3cm]{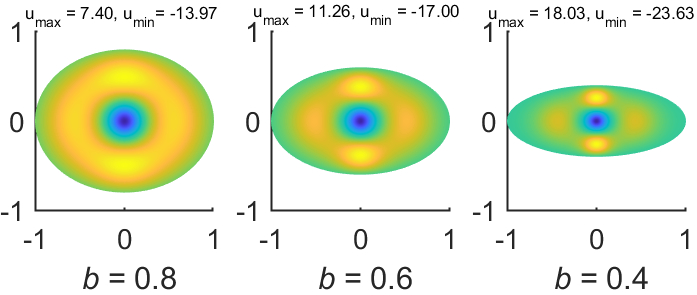}}\\ \subfigure[VII]{\includegraphics[width=7.15cm,height=3cm]{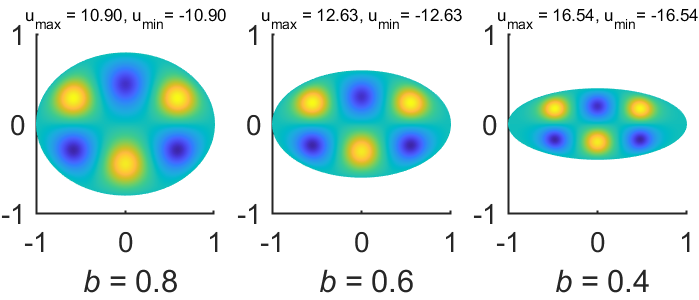}} \subfigure[VIII]{\includegraphics[width=7.15cm,height=3cm]{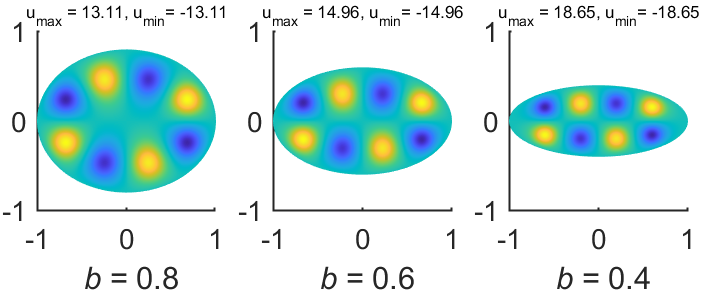}}\\ \subfigure[IX]{\includegraphics[width=7.1cm,height=3cm]{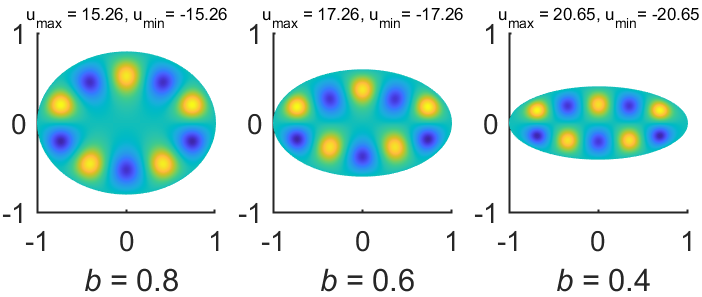}}
\caption {Multiple solutions of \eqref{eq.4.1} with varying $b$.}\label{example21}
\end{center}
\end{figure}

\clearpage
According to \cite{xie2004}, \eqref{eq.4.1} has infinitely many solutions. Moreover, rotating a solution $u$ by integer times of 90 degrees around the center, the resulting function rot($u$) is still a solution of \eqref{eq.4.1}; If $u$ is a solution of \eqref{eq.4.1}, so is $-u$.
{\color{black}Fig.\ref{example21type} presents some multiple solutions calculated using the proposed method, where $M = N = 25$ in \eqref{20240112eq2.1} is chosen based on the Legendre-Fourier scheme described in Section \ref{sect2}.
These solutions agree well with the results given in \cite{2022Two}}. $J(u)$ and $J^\star(u):=J(u)/J(u_I)$ corresponding to these solutions are plotted in Fig.~\ref{20241014example22}. 
These solutions with varying $b$ are also shown in Fig.~\ref{example21}.\vspace{-5cm}

{\color{black} Theoretically, for a given solution $u$, having more oscillations in the short axis direction leads to larger 
 $|\nabla u|$ (thus larger energy) than that having the same oscillations in the long axis direction. The type-II solution shown in Fig.\ref{example21type}(b), and the solution obtained by rotating it 90 degrees, have the same energy. However, using these two different solutions as initial values to calculate solutions of the case $b<1$, one obtains solutions with different energies. The results of relative energy $J^\star(u):=J(u)/J(u_I)$ for the two cases that use solutions obtained by rotating type-II and type-III 90 degrees as initial guesses for cases $b<1$ are presented in the right plot (the blue curve and the red curve) of Fig..~\ref{20241014example22}, which justify the aforementioned theoretical understanding.\vspace{-8.5cm}
}
\newpage
{\em \textbf{Example 3.}} Motivated by \cite{xie100810411}, we consider the following singularly perturbed Ginzburg--Landau problem
\begin{equation}\label{2024082903}
\begin{cases}
\delta \triangle u - u + u^{3} = 0 \quad\quad\quad \textrm{in}\;\Omega,\\
u = 0   \quad\quad\quad\quad\quad \quad\quad\quad\;\,  \textrm{on}\;\partial\Omega,
\end{cases}
\end{equation}
where $0 < \delta = \varepsilon^2 \ll 1$. The variational functional corresponding to \eqref{2024082903} is
\begin{equation}\label{eq33}
J(u) = \int_{\Omega}(\frac{\delta}{2}|\nabla u|^2 + \frac{1}{2}u^2-  \frac{1}{4}u^4)dx,  \quad\quad\quad  u \in H^1_0(\Omega).
\end{equation}

Some theoretical results presented in \cite{GUI19991,Cao1999OnTE,GUI200047,Gui2000OnMM} and numerical results obtained using newly proposed method in this paper are given as follows:
\begin{figure}[H]
\begin{centering} \subfigure[I]{\includegraphics[width=3.90cm,height=6.4cm]{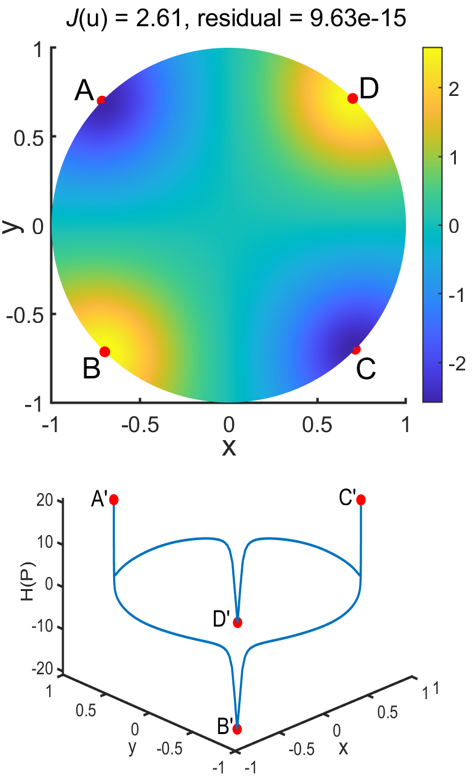}}\;\; \subfigure[II]{\includegraphics[width=3.93cm,height=6.4cm]{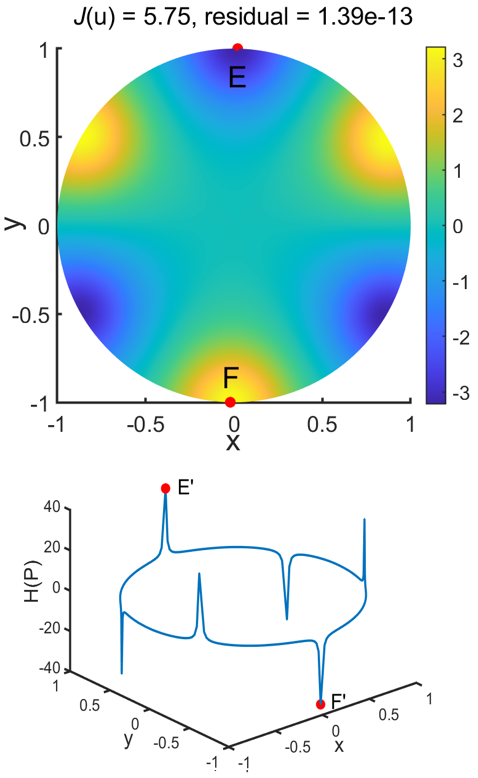}}\;\; \subfigure[III]{\includegraphics[width=3.95cm,height=6.4cm]{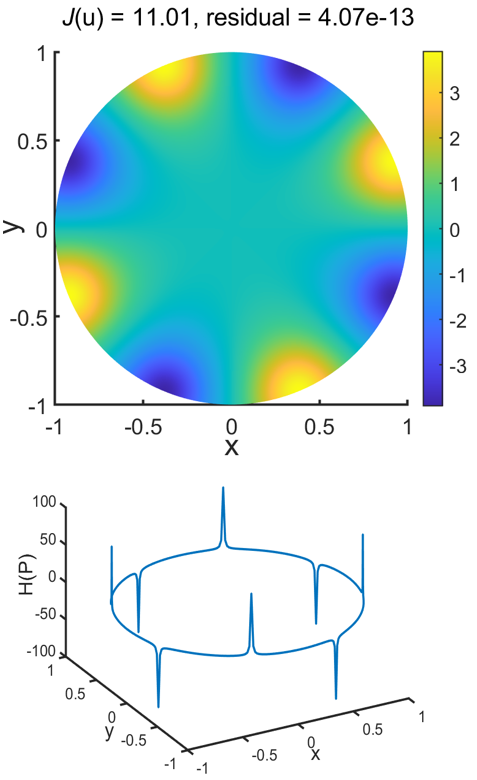}}\;\;
\subfigure[IV]{\includegraphics[width=3.85cm,height=6.3cm]{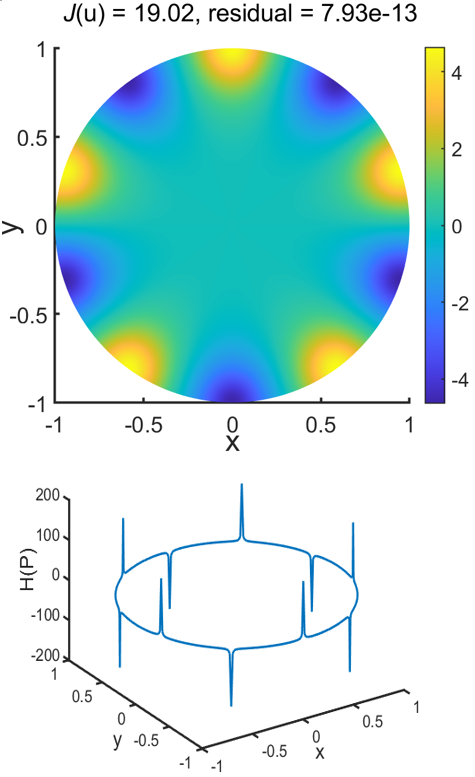}}\;\;
\subfigure[X]{\includegraphics[width=3.9cm,height=6.3cm]{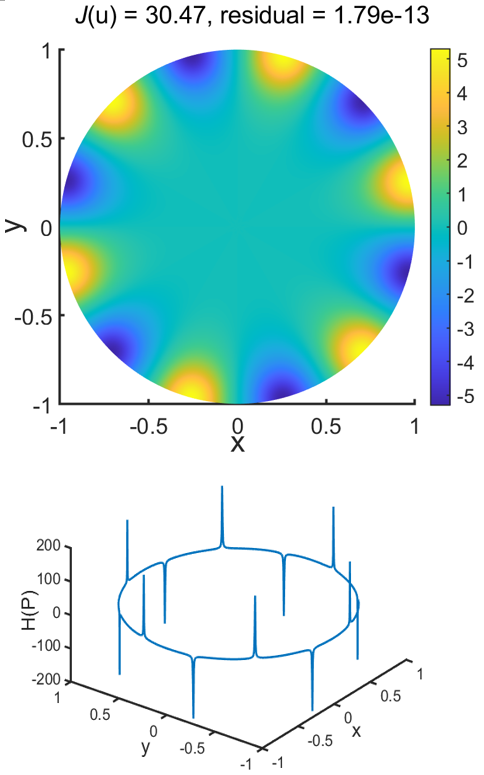}}\;\;
\subfigure[VI]{\includegraphics[width=3.9cm,height=6.3cm]{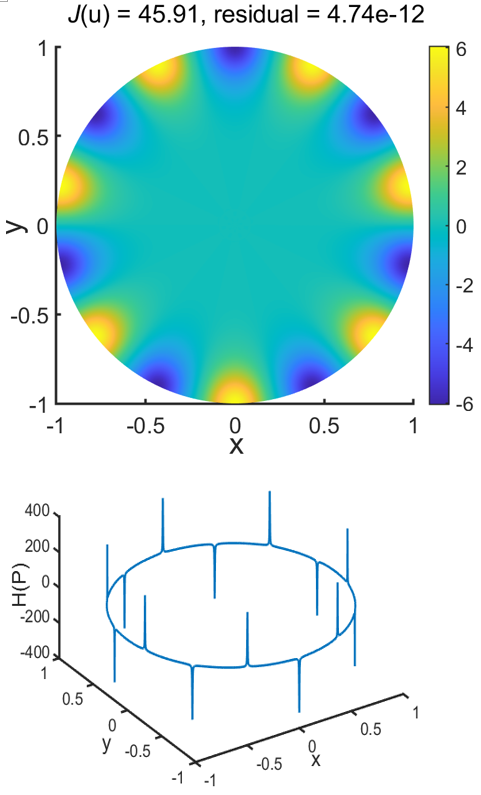}}\;\;
\caption{Multiple solutions to \eqref{2024082903} and their mean curvatures $H(P)$ at $\partial \Omega$ with $\delta = 0.2$. For example, for the type-I solution, the peaks ($A,B,C$ and $D$) at $\partial \Omega$ are located at the maximum or minimum of $H(P)$ ($A',B',C'$ and $D'$).}\label{example31}
\end{centering}
\end{figure}
\begin{itemize}
\item There exist some solutions with some peaks on $\partial\Omega$, and these peaks are located at which the mean curvature of the solution surface on $\partial\Omega$ is the maximum or minimum. These phenomena can be seen in Fig.\ref{example31}. In Figs.\ref{2024101501example34}-\ref{2024101502example34}, the type-I and type-II solutions are shown with varying $b$, where these peaks go to endpoints with decreasing $b$. In addition, $J(u)$ corresponding to multiple solutions presented in Fig.\ref{example31} are plotted with varying $b$ (see Fig.\ref{example32}).

\begin{figure}[H]
\begin{center}	
\subfigure[$b = 0.8$]{\includegraphics[width=3.4cm,height=3.75cm]{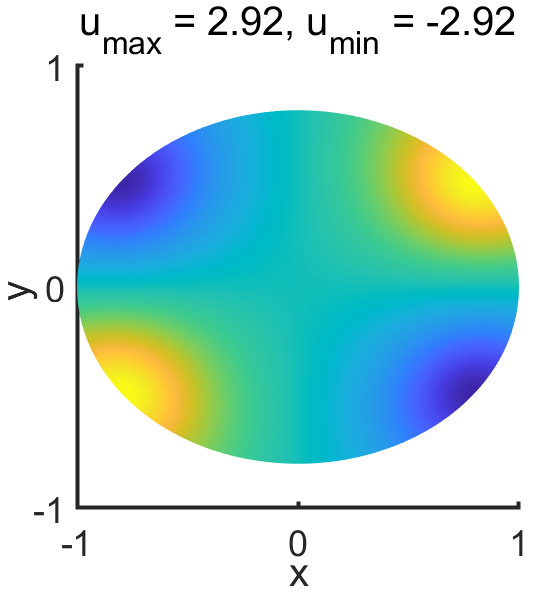}}\quad\quad \subfigure[$b = 0.6$]{\includegraphics[width=3.4cm,height=3.7cm]{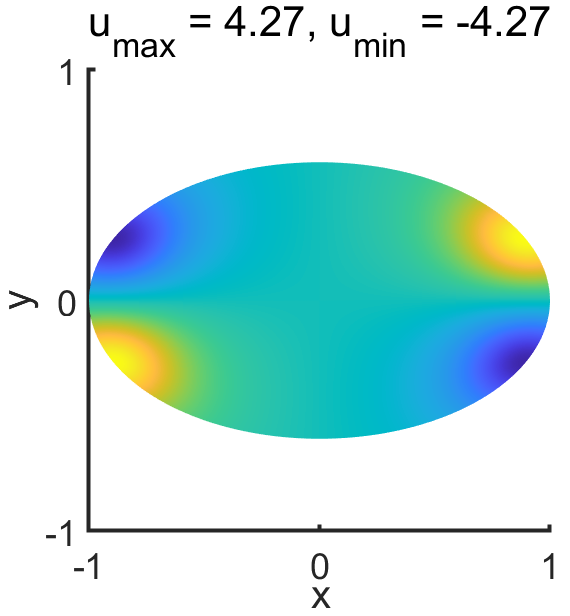}}\quad\quad \subfigure[$b = 0.4$]{\includegraphics[width=3.4cm,height=3.75cm]{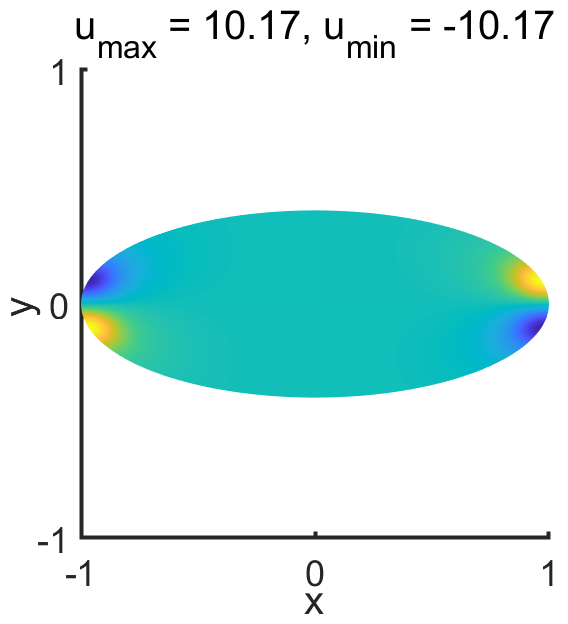}}
\caption{The type-I solutions to \eqref{2024082903} with varying $b$.}\label{2024101501example34}
\subfigure[$b = 0.8$]{\includegraphics[width=3.4cm,height=3.75cm]{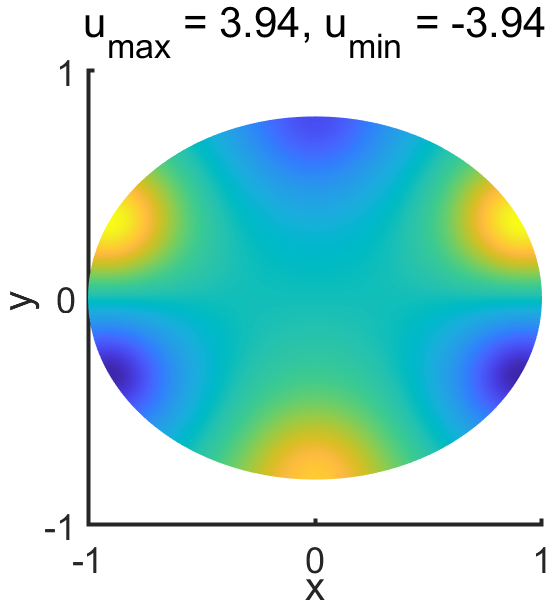}}\quad\quad \subfigure[$b = 0.6$]{\includegraphics[width=3.4cm,height=3.75cm]{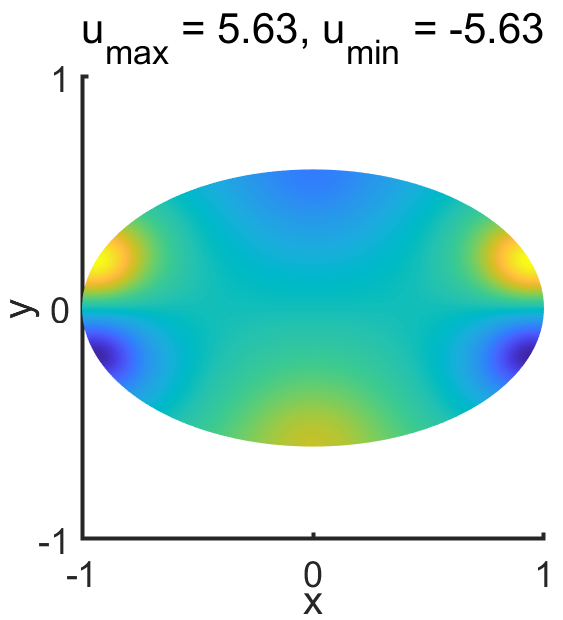}}\quad\quad \subfigure[$b = 0.4$]{\includegraphics[width=3.4cm,height=3.75cm]{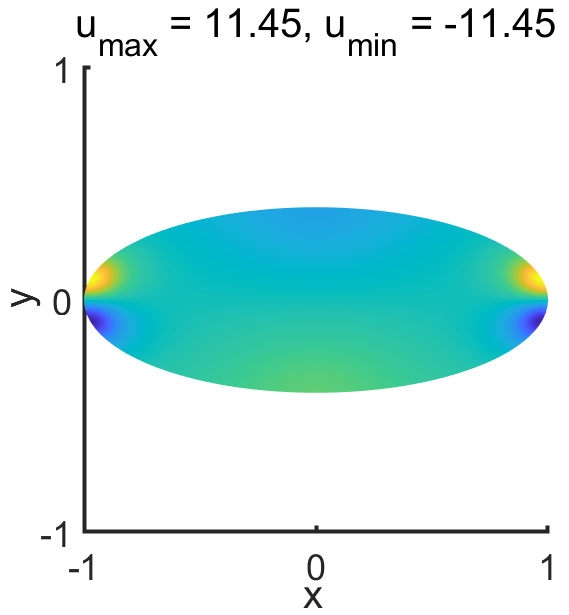}}
\caption{The type-II solutions to \eqref{2024082903} with varying $b$.}\label{2024101502example34}
\end{center}
\end{figure}

\begin{figure}[!h]
\vspace{-0.5cm}
\begin{centering}
\includegraphics[width=5.1cm,height=4.8cm]{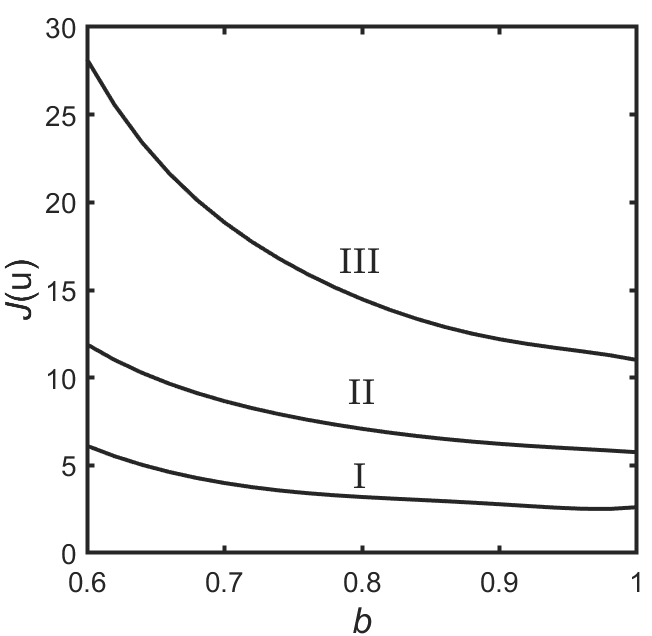}\quad\quad\quad
\includegraphics[width=5.1cm,height=4.8cm]{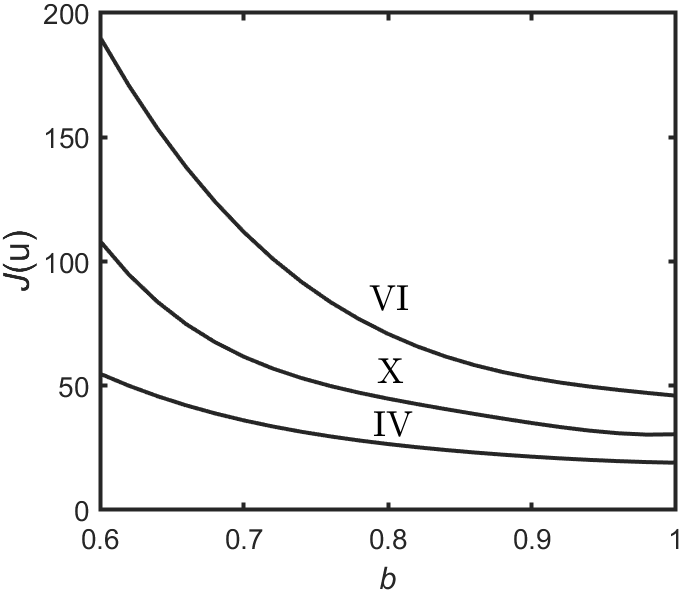}
\caption{$J(u)$ vs $b$ with $\delta = 0.2$.}\label{example32}
\vspace{-0.5cm}
\end{centering}
\end{figure}

\item For $\delta$ sufficiently small, the solution $u$ has at most one local maximum and it is achieved at exactly one point $P_{\delta}$ in $\Omega$. Moreover, $u \rightarrow 0$ in $\Omega \setminus P_{\delta}$. In other words, the peak of the solution becomes sharper as $\delta \rightarrow 0$ (i.e. point-condensation phenomenon), which can be seen in Fig.\ref{example33}.

\begin{figure}[H]
\begin{centering}
\subfigure[$\delta = 2\times10^{-2}$]{\includegraphics[width=4.3cm,height=3.9cm]{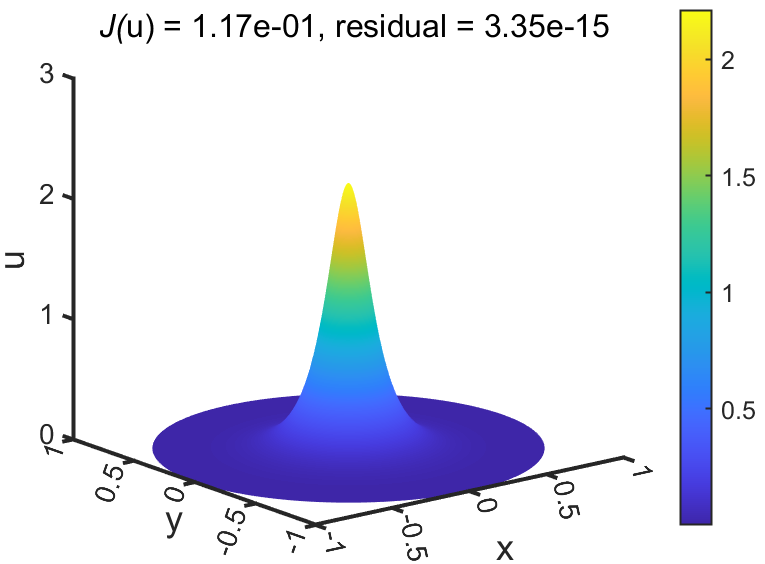}}\;\;\;\;
\subfigure[$\delta = 2\times 10^{-4}$]{\includegraphics[width=4.3cm,height=3.9cm]{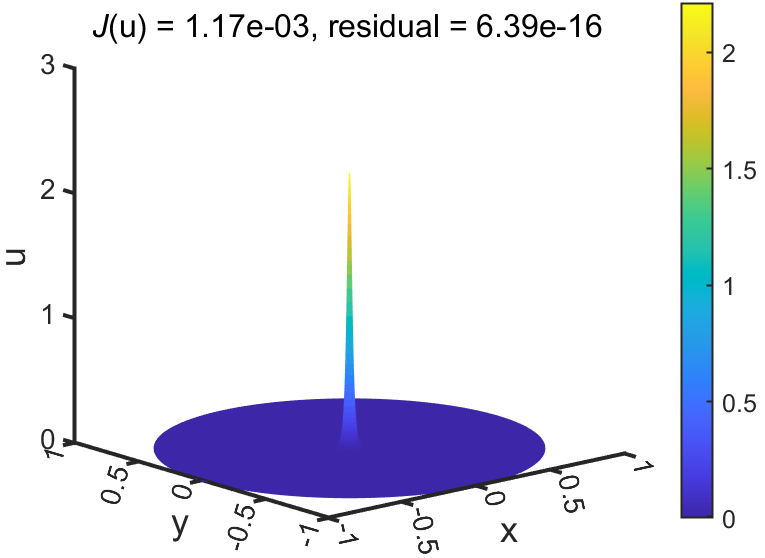}}\;\;\;\;
\subfigure[$\delta = 2\times10^{-6}$]{\includegraphics[width=4.3cm,height=3.9cm]{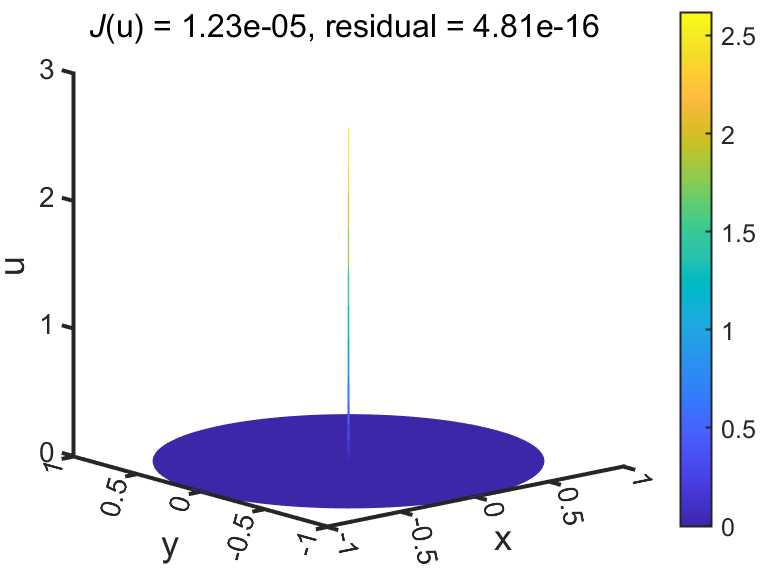}}
\caption{The single-peak solutions to \eqref{2024082903} with different values of $\delta$.}\label{example33}	
\end{centering}
\end{figure}

\item When $\delta \rightarrow 0$, there exist solutions with some peaks located at $\partial\Omega$ or the interior region, which can be seen in Fig.\ref{example34}.
\end{itemize}
\begin{figure}[H]
\begin{centering}	
\subfigure[I]{\includegraphics[width=4.4cm,height=4.1cm]{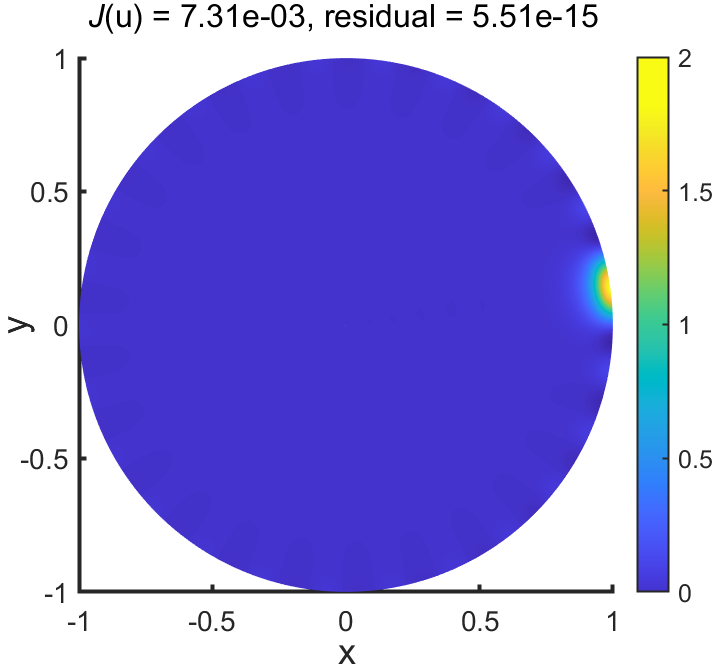}}\;\;\;	
\subfigure[II]{\includegraphics[width=4.4cm,height=4.1cm]{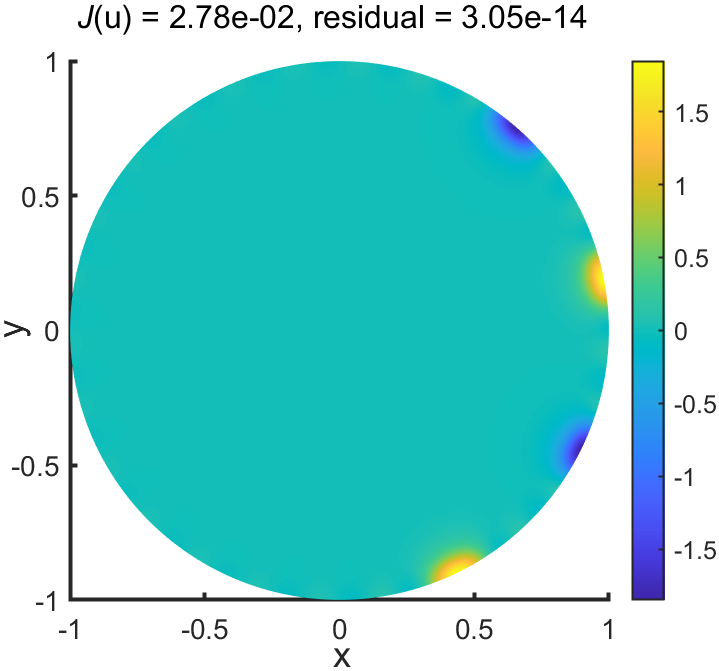}}\;\;\;			
\subfigure[III]{\includegraphics[width=4.4cm,height=4.1cm]{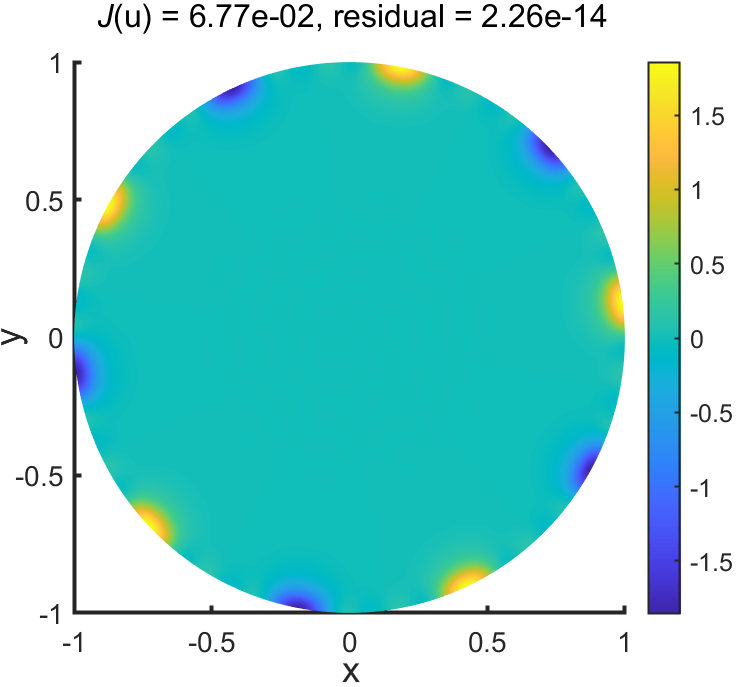}}\;\;\\ 
\subfigure[IV]{\includegraphics[width=4.4cm,height=4.1cm]{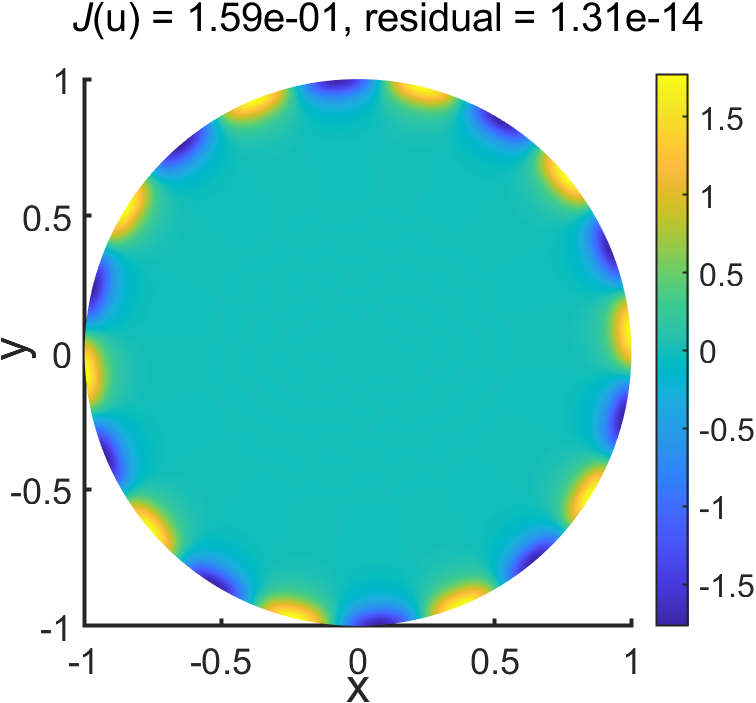}}\;\;\; 
\subfigure[V]{\includegraphics[width=4.4cm,height=4.1cm]{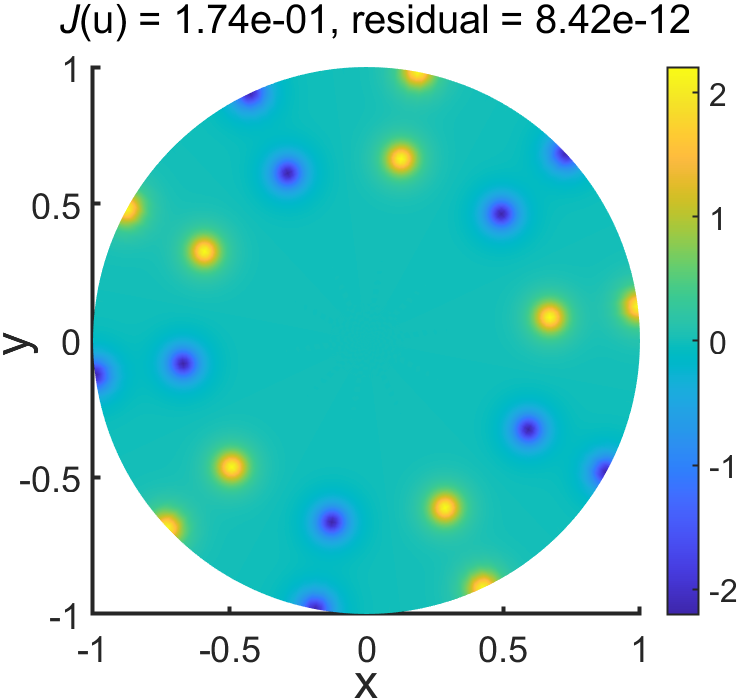}}\;\;\; 
\subfigure[VI]{\includegraphics[width=4.4cm,height=4.1cm]{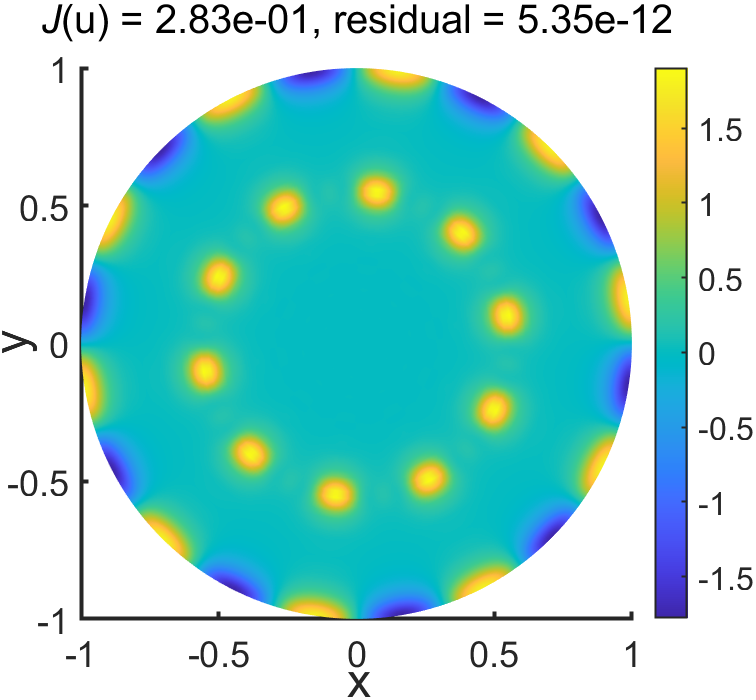}}\;\;\\

\caption{Multiple solutions to \eqref{2024082903} with $\delta = 2\times 10^{-3}$.}\label{example34}
\end{centering}
\end{figure}

\section{Conclusions}\label{sect5}
In this paper, an improved adaptive orthogonal basis deflation method for computing multiple solutions of nonlinear elliptic equations is proposed, {\color{black}through which the impact of domain changes on multiple solutions of several nonlinear elliptic equations are numerically investigated}. The improved adaptive orthogonal basis method accurately validates some theoretical results on multiple solutions of these examples. This not only illustrates the efficiency of the proposed method but also opens new avenues for integrating numerical computation into theoretical analysis.

 Although only nonlinear elliptic partial differential equations are considered in this study, our proposed numerical method {\color{black}can apply to} 
 more general non-convex partial differential equations arising in diverse and practical problems, e.g., material science problems~\cite{yu_kinetichydrodynamic_2007, e_minimum_2004, zhangDynamicTransitionsLandauBrazovskii2014, wang_energystable_2019, yinConstructionPathwayMap2020}, fluid dynamics problems~\cite{wan_model_2015, wan_dynamics_2017} and topology optimization and electronic design problems~\cite{wangEfficientUnconditionallyStable2022, liProvablyEfficientMonotonicdecreasing2022, yue_efficient_2024}. In addition, it provides a useful numerical tool to verify more theoretical results presented in published literature, e.g., \cite{Wei2023,WeiWu2023,LIU2022509}. {\color{black}We will continue to apply and extend the proposed multiple-solution approach to address more challenging problems in the future}.

\section*{Acknowledgements}
The work of H. Yu was partially supported by National Key R\&D Program of China (grant number 2021YFA1003601)
and the National Natural Science Foundation of China (grant number 12171467, 12161141017). 


\end{document}